  \documentclass[final,3p, review,11pt]{elsarticle}
  \usepackage{pifont}

\usepackage{lineno,hyperref}
\modulolinenumbers[5]

\usepackage{amssymb}
\usepackage{amsfonts}
\usepackage{amsmath,amsthm}
\usepackage{makeidx,pstricks,latexsym,bm,graphicx}
\usepackage{subfigure}
\usepackage{mathrsfs}
\usepackage{multirow}
 \usepackage{txfonts}

\begin{document}
\begin{frontmatter}

\title{An improved isogeometric analysis method for trimmed geometries}

\author[A]{Jinlan Xu} \ead{jlxu@hdu.edu.cn}

\author[A]{Ningning Sun}

\author[A]{Laixin Shu}

\author[label4]{Timon Rabczuk}

\author[A]{Gang Xu\corref{cor} } \ead{gxu@hdu.edu.cn}

\cortext[cor]{Corresponding author}
\address[A]{School of Computer Science and Technology, Hangzhou Dianzi
  University, Hangzhou 310018, P.R.China}
  \address[label4]{Institute of Structural Mechanics, Bauhaus-University Weimar, Marienstr. 15, D-99423 Weimar, Germany}

\begin{abstract}

Trimming techniques are efficient ways to generate complex geometries in Computer-Aided Design(CAD). In this paper, an improved isogeometric analysis(IGA) method for trimmed geometries is proposed. We will show that the proposed method reduces the numerical error of physical solution by 50\% for simple trimmed geometries, and the condition number of stiffness matrix is also decreased. Furthermore, the number of integration elements and integration points involved in the solving process can be significantly reduced compared to previous approaches, drastically improving the computational efficiency for IGA problems on the trimmed geometry. Several examples are illustrated to show the effectiveness of the proposed approach.

\end{abstract}


\end{frontmatter}

\section{Introduction}
\label{intro}

Isogeometric analysis(IGA)~\cite{Hughes-Cottrell-2005}, was motivated by unifying Computer Aided Design (CAD) and Finite Element Analysis (FEA). Geometries in CAD are usually represented by splines(B-splines, NURBS, T-splines, PHT-splines for instance), with the geometries in FEA are commonly based on Lagrange polynomials. These two different geometry descriptions introduce inconsistencies in CAD and CAE designs which requires reapproximating the CAD geometries in CAE. This does not only introduce errors in the geometry but increase the entire design-to-analysis time. It was demonstrated in~\cite{AnalysisAware-2010,xu:gmp10,xu:cmame2011} that the mesh quality has a big impact on the analysis results, and meshing operation occupies a large percentage in the entire analysis procedure. Isogeometric analysis unifies the geometry representation of design and analysis, by using the same CAD spline functions in CAE simplifying the design-analyze process, and ensuring the exact geometry during the analysis. If a high precision numerical solution is requested, mesh refinement is inevitable. In FEA, posterior error is often used to guide the refinement, and the refinement based on the mesh is sometimes not appropriate, so re-meshing will be needed which have to be interact with original model. In practical engineering analysis, this is a severe bottleneck. IGA solved these problems by applying same bases to represent the geometry and the physical phase, which
results in a direct design-to-analysis process without the intermediate step of mesh generation.
The geometry is exact at the coarsest level, thus eliminating the geometrical errors, and can be refined by knot insertion
or order elevation. Refinement at any level can take place completely within the analysis framework, which eliminates the necessity to communicate with the geometry.

Most of CAD models can not be represented by a single tensor-product spline surface but
several patches of spline surfaces are needed ~\cite{xu:cmame2015}. However, if a geometry is complex and is required to be smooth in the interior, it is not easy to construct such a geometry with multiple patches of spline surfaces. In such case, trimming techniques are usually employed. Other ways are to approximate the geometry with T-splines~\cite{Tpline-2003,Planar-Tspline-2014}, PHT-splines~\cite{PHT-2008,Chan-Timon-2017}, THB-splines~\cite{Falini-Bert-2015}, or LR B-splines~\cite{Johannessen-Dokken-2014}. If trimming techniques are used, normal elements and trimmed elements will be considered separately during isogeometric analysis. Kim~\cite{Hyun-Jung Kim_2009} proposed a method to solve this problem. Schmidt et al.\cite{Schmidt-2012} proposed a reconstruction method using geometric bases to evaluate the finite element constituents of trimmed elements. This method covers both bases of a single patch and multi-patches. Shen~\cite{Shen-2014} introduced a method to convert trimmed NURBS surfaces to subdivision surfaces, and their method can produce gap-free models which is mandatory for numerical analysis. Moreover, the resulting models are $G^1$ continuous between two adjacent surfaces.  Zhu et al.~\cite{Zhu-Hu-2016} proposed a spline called B++ spline, to express the trimmed NURBS patch in an analytic form. They solved the problem of implementing essential boundary conditions in isogeometric analysis. The basis functions of B++ spline satisfy the Kronecker delta property which allow imposing essential boundary condition strongly, and this is similar with FEM. Other interesting approaches on isogeometric analysis for trimmed surfaces, can be found in~\cite{Kang-2015,Michael-2013,Beer-2014,Martin-2014,Wang-2015,ZhuXF-2016,Marussig-2016} and references therein. In this paper, we improve the method proposed by Kim~\cite{Hyun-Jung Kim_2010} which is among the most efficient approaches for trimmed geometries in our opinion.

The paper is organized as follows. In section~\ref{pre}, we summarize the basics of the IGA formulation on trimmed geometries presented in~\cite{Hyun-Jung Kim_2010}. In section~\ref{integration-pts}, we describe our method to deal with trimmed element in details. Section~\ref{example} gives several examples of our proposed method, in comparison to the method in~\cite{Hyun-Jung Kim_2010} are also presented. We end this paper with conclusions in section~\ref{conlude}.

\section{Preliminaries}
\label{pre}

NURBS are the most common basis functions for representing free-form objects. However, tensor product form of NURBS makes the gap-free(and overlap-free) representation of more
complex objects cumbersome. Trimming techniques eliminate this limitation of NURBS. Since the isogeometric analysis is based on the tensor-product mesh structure,
geometries with trimming curves have to be decomposed into several NURBS surfaces in traditional isogeometric analysis.  Kim et al. proposed an IGA formulation on trimmed geometries directly~\cite{Hyun-Jung Kim_2009,Hyun-Jung Kim_2010}, which is summarized subsequently.

\subsection{Trimmed geometry}
\label{trim-geo}

\begin{figure}[!htb]
\small
\begin{center}
  \begin{tabular}{cc}
  \includegraphics[width=2.50in]{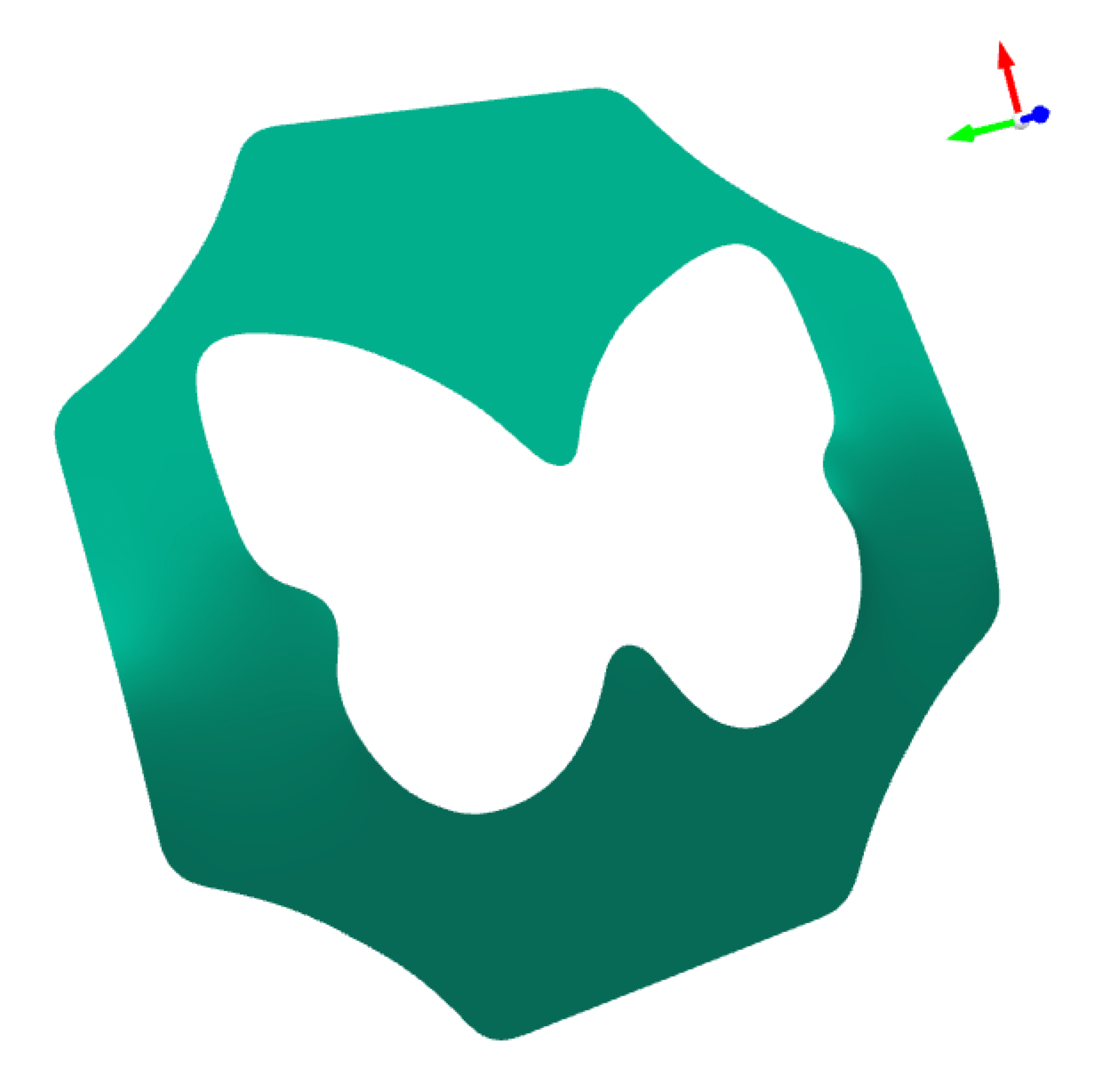} &
  \includegraphics[width=2.50in]{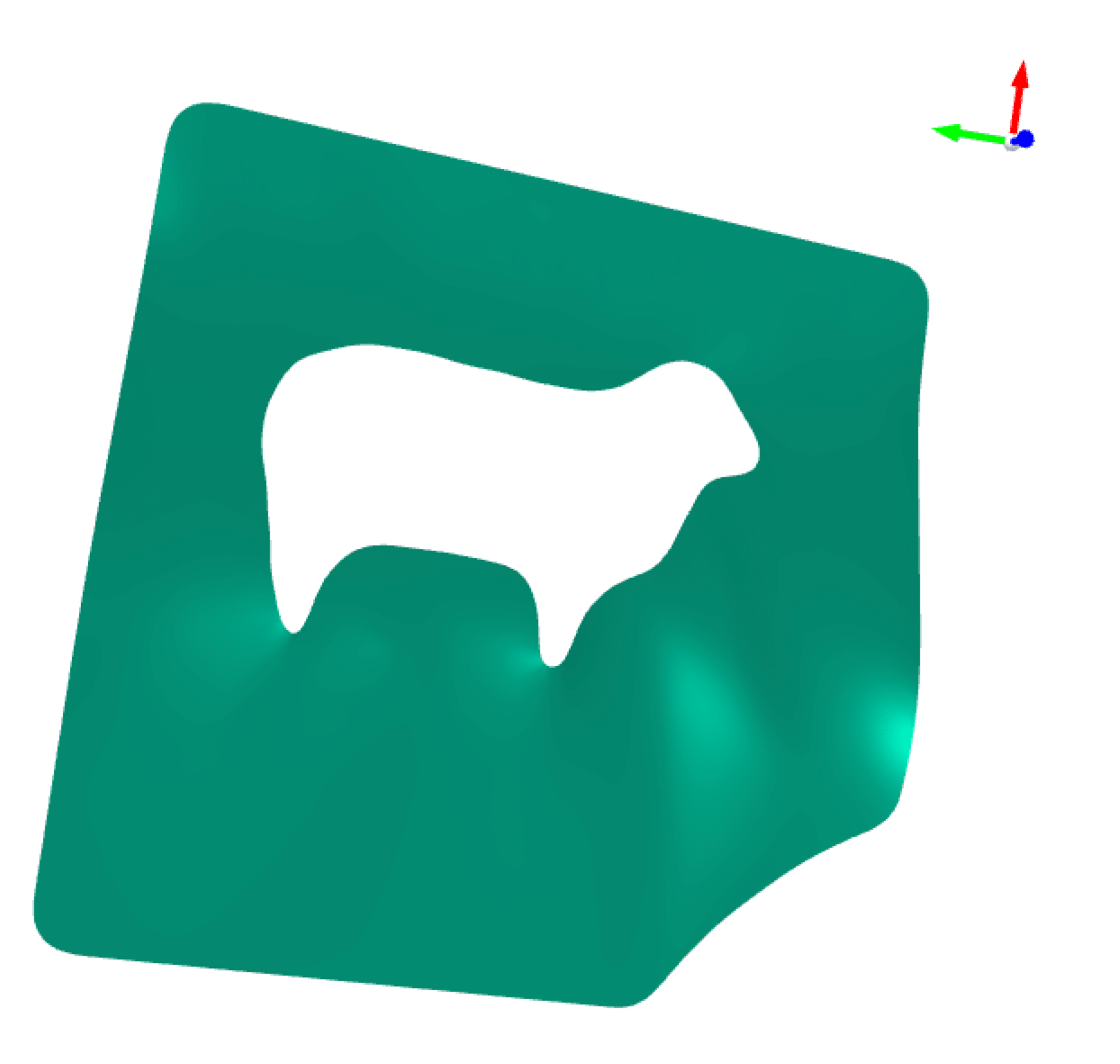} \\
  (a) & (b)
  \end{tabular}
\caption{\label{fig-TrimSurf} Two examples of trimmed NURBS surfaces.}
\end{center}
\end{figure}

In CAD, geometries are usually represented by NURBS. However, it is not trivial to represent the geometry exactly with one NURBS surface or NURBS solid for complex geometries. Trimming techniques employ NURBS curves to trim unwanted parts of geometries as shown in Fig.\ref{fig-TrimSurf}. The resulting geometry is called a trimmed geometry, but there is no mathematical relation between the trimming curve and the NURBS surface. For a trimmed surface, the CAD files contain the surface information in the  parametric space and physical space. We only consider the parametric data, and the operation on the parametric data will be mapped onto the geometric data.
The trimming technique not only simplifies the construction of complex models, but also keep the smoothness of the untrimmed parts.

\subsection{Isogeometric analysis on trimmed geometry}
\label{trim-IGA}

As there are many trimmed patches in the CAD models, trimmed geometries is a key problem for IGA.  Kim et al. \cite{Hyun-Jung Kim_2009} proposed the first approach to analyze the geometry directly on trimmed surfaces.

 \begin{figure}[!htb]
\small
\begin{center}
  \begin{tabular}{ccc}
  \includegraphics[width=1.90in]{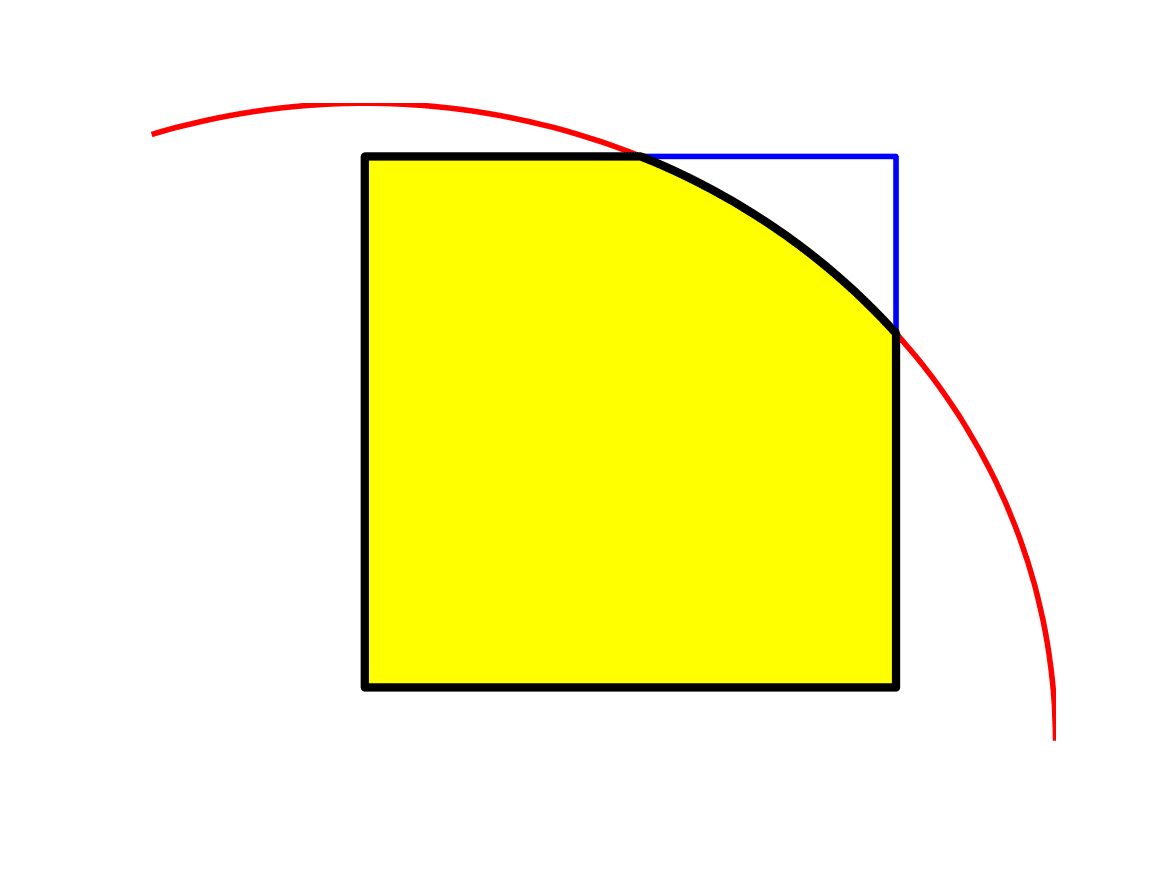} &
  \includegraphics[width=1.90in]{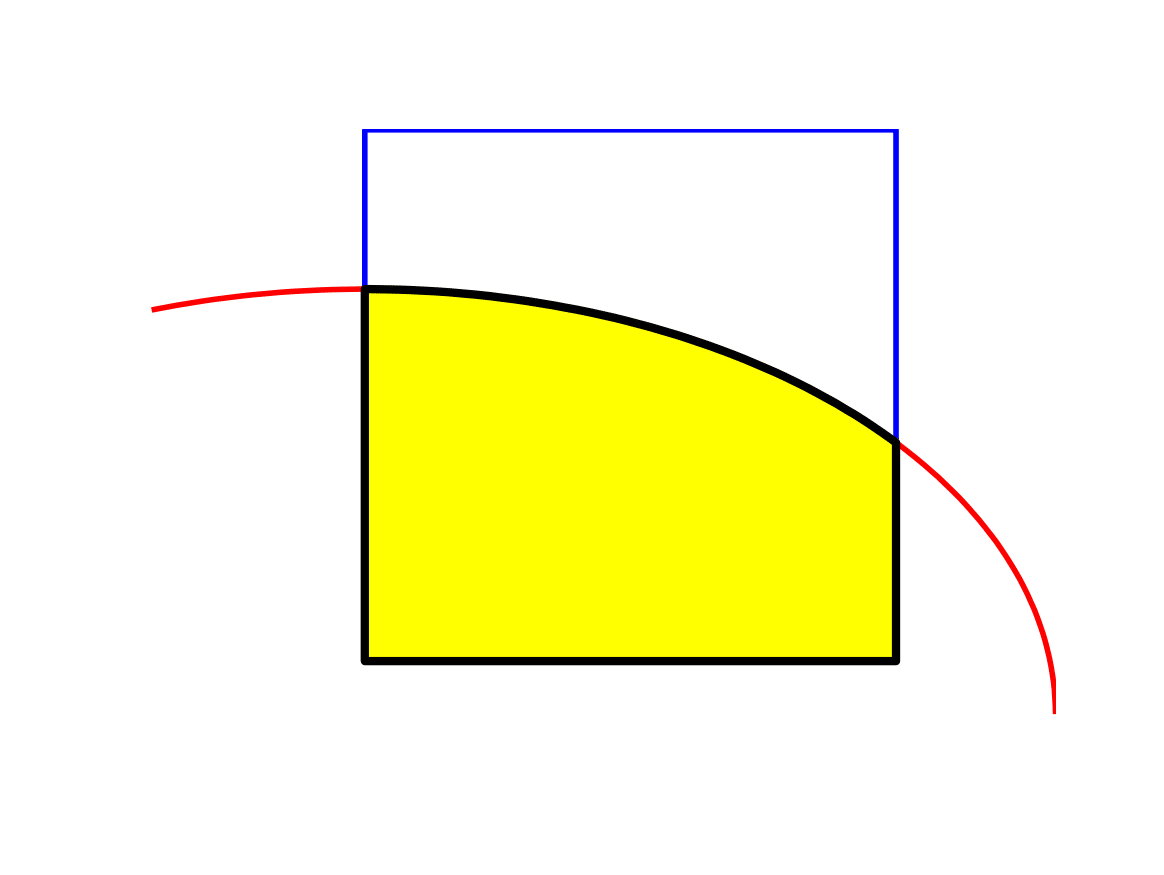} &
  \includegraphics[width=1.66in]{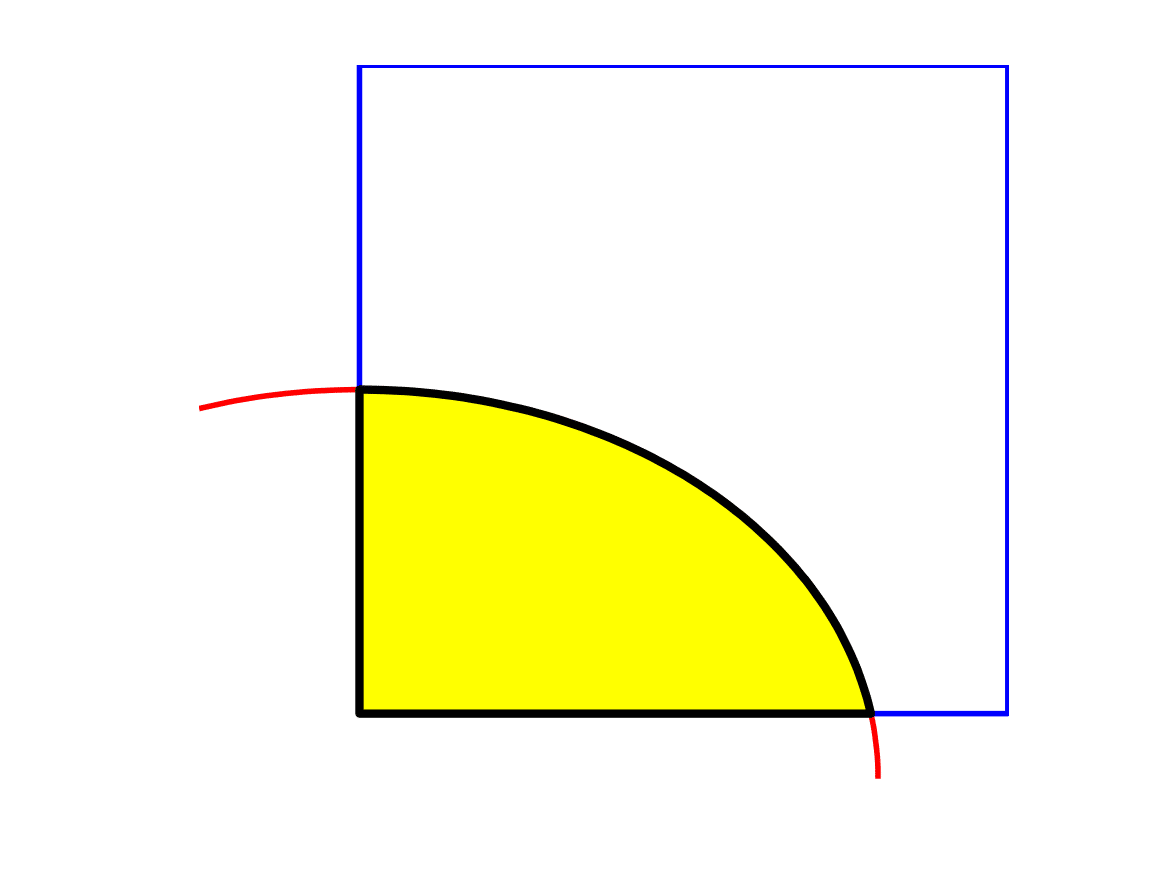} \\
  (a) & (b) & (c)
  \end{tabular}
\caption{\label{fig-eletype} Three types of elements:(a)type A;(b)type B;(c)type C.}
\end{center}
\end{figure}

\begin{figure}[!htb]
\small
\begin{center}
  \begin{tabular}{cc}
  \includegraphics[width=1.40in,clip]{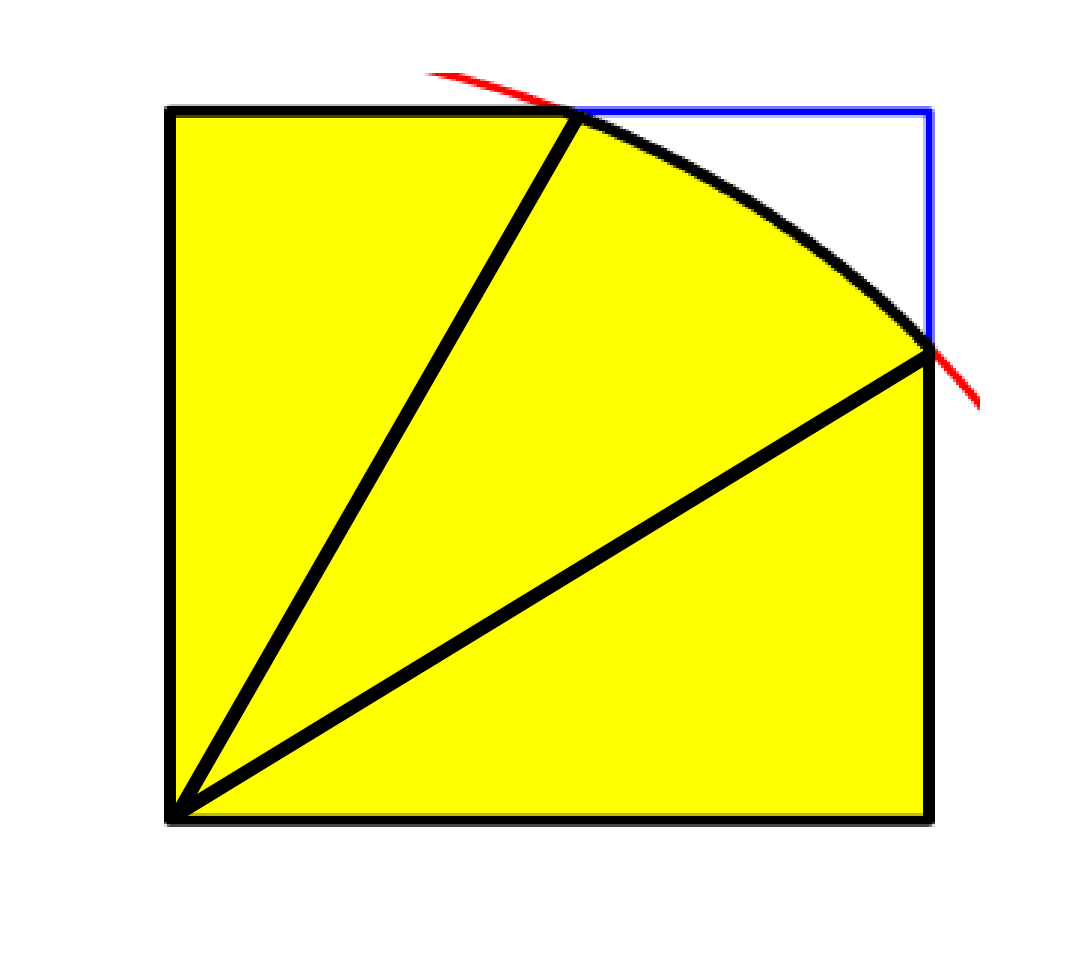} &
  \includegraphics[width=1.4in,clip]{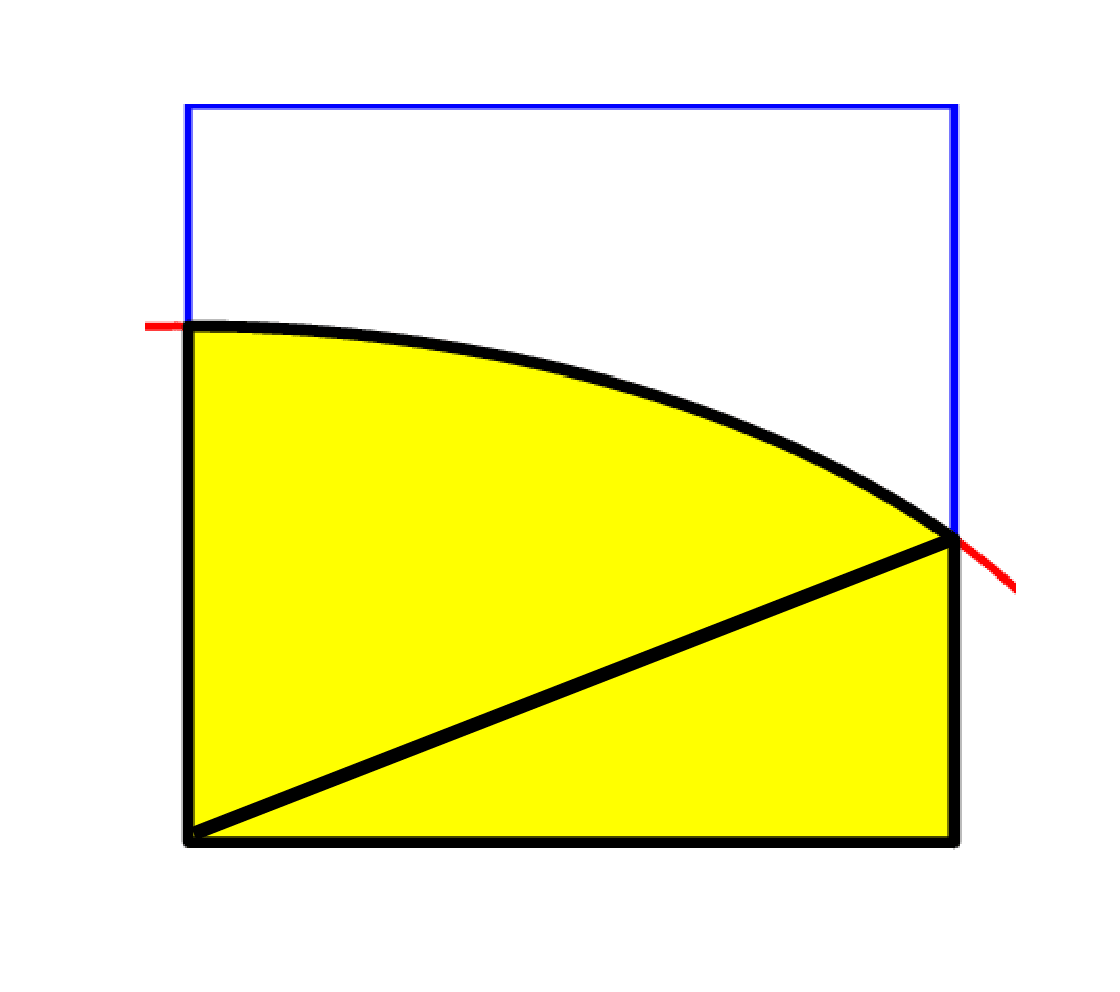} \\
  (a)type A & (b)type B
  \end{tabular}
\caption{\label{fig-seg-eletype} Segmentation of elements with type A and type B.}
\end{center}
\end{figure}

Therefore, trimmed elements are divided into three types, as depicted in Fig.~\ref{fig-eletype}.
\begin{itemize}
\item Type A is the curved pentagon;
\item Type B is the curved quadrilateral;
\item Type C is the curved triangle;
\end{itemize}
where the curved edge is part of the trimming curve. For each type, they segment the trimmed element to triangular ones. The triangular elements are further subdivided to two cases, i.e. a(normal) triangle with straight edges and a triangle with two straight and one curved edge, see Fig.~\ref{fig-seg-eletype}. Let us denote these triangles as $T$ and $\tilde{T}$ respectively. For the 'normal' triangle $T$ the integration points for the linear triangle(triangular element) are chosen. If $T$ has a curved edge, a series of triangulation is performed on $T$ in order to map $T$ to a rectangle.

\subsubsection{Imposition of essential boundary condition}

In isogeometric analysis, essential boundary conditions can not be imposed as in FEM, because NURBS basis functions do not satisfy the Kronecker delta property. For homogeneous essential boundary conditions, the coefficients of basis functions corresponding to boundary are set to zero. The imposition of non-homogeneous essential boundary conditions require special techniques such as modification of the weak form or the solution of an interpolation problem at the boundary; see e.g.~\cite{Martin-2014}.

In trimmed isogeometric analysis, additional challenges occur for imposing essential boundary condition. The boundary conditions need to be imposed on the trimming curves but the degree of freedom(DOF) is defined on the NURBS surface. Furthermore there is no mathematical relationship between these two representations. In~\cite{Hyun-Jung Kim_2010}, they use Lagrange multiplier method to impose essential boundary conditions on trimming curves. And we use the same method which is described in Section~\ref{example}.

\section{Improved isogeometric analysis on trimmed geometry}
\label{integration-pts}

In this paper, we propose an improved method for IGA on trimmed geometries, by which improves the efficiency and accuracy of trimmed IGA.

\begin{figure}[!htb]
\small
\begin{center}
 \includegraphics[width=2in,clip]{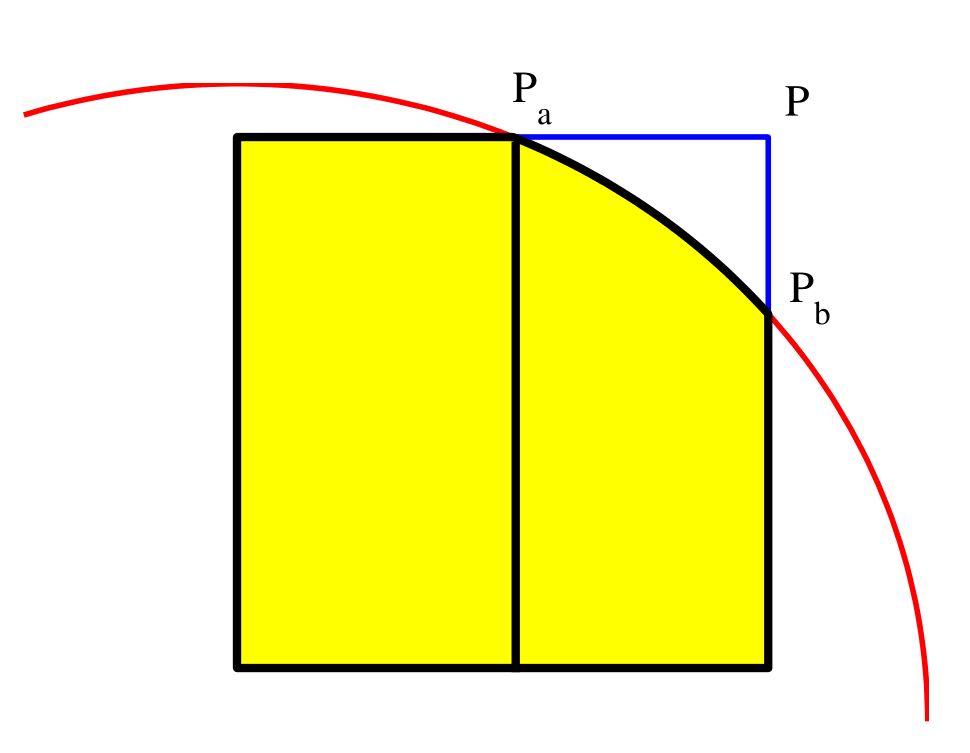}
\caption{\label{fig-typeA-decompose} The trimmed element of type A is decomposed to two quadrilaterals. }
\end{center}
\end{figure}

The main contribution of our method is to modify the integration rules for the trimmed elements. For type C in ~\cite{Hyun-Jung Kim_2010}, a similar method is applied to generate integral points on the curved triangle, but for type B~\cite{Hyun-Jung Kim_2010}, a mapping from a rectangle to the curved quadrilateral element is used which avoids the triangular decomposition of the curved quadrilateral element. For type A, a decomposition is adopted but it is different from the method in ~\cite{Hyun-Jung Kim_2010}. The curved pentagon is segmented to two quadrilaterals, one with a curved edge and the other is rectangular. Fig.~\ref{fig-typeA-decompose} shows the decomposition of type A in our method.

The segmentation of trimmed elements of type A can be chosen on the basis of the intersection points. Suppose $P_a$ and $P_b$ are two intersection points, where $P_b$ is closer to the corner point which is trimmed out, then the trimmed element is segmented at point $P_a$.

\begin{figure}[!htb]
\small
\begin{center}
  \begin{tabular}{cc}
  \includegraphics[width=1.78in,clip]{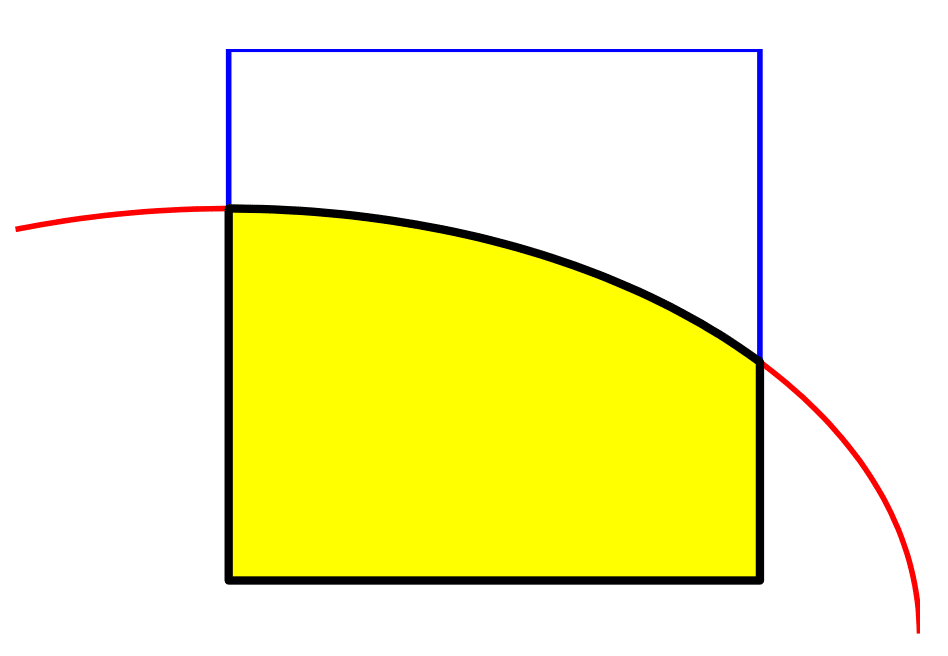} &
  \includegraphics[width=1.78in,clip]{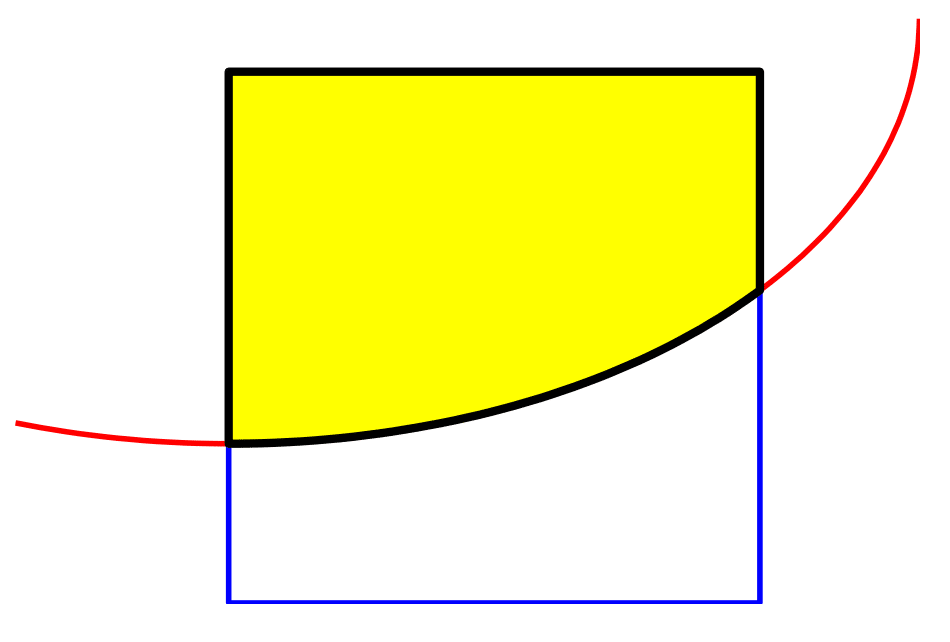} \\
  (a) & (b) \\
  \includegraphics[width=1.48in,clip]{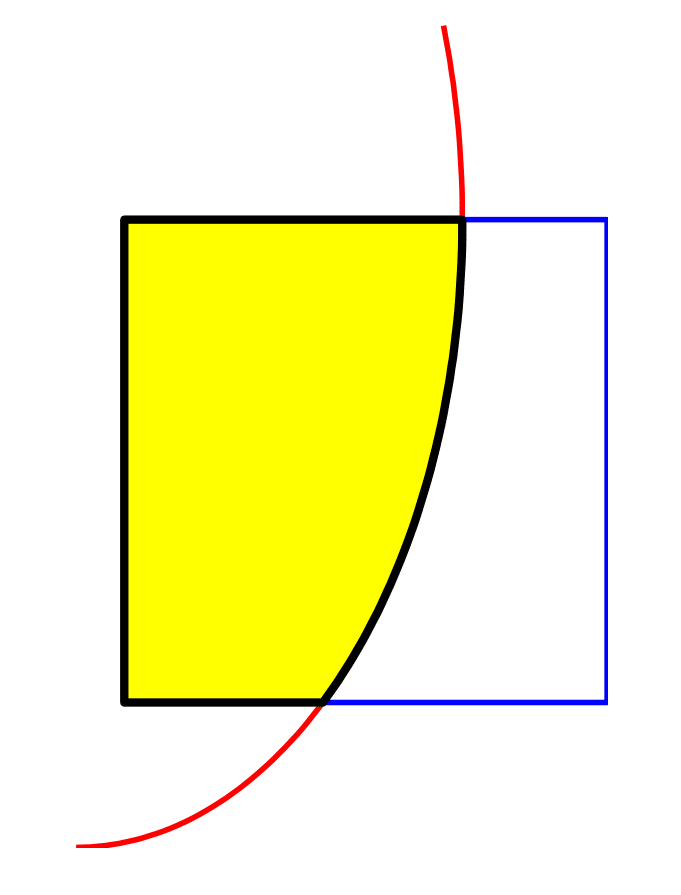} &
  \includegraphics[width=1.30in,clip]{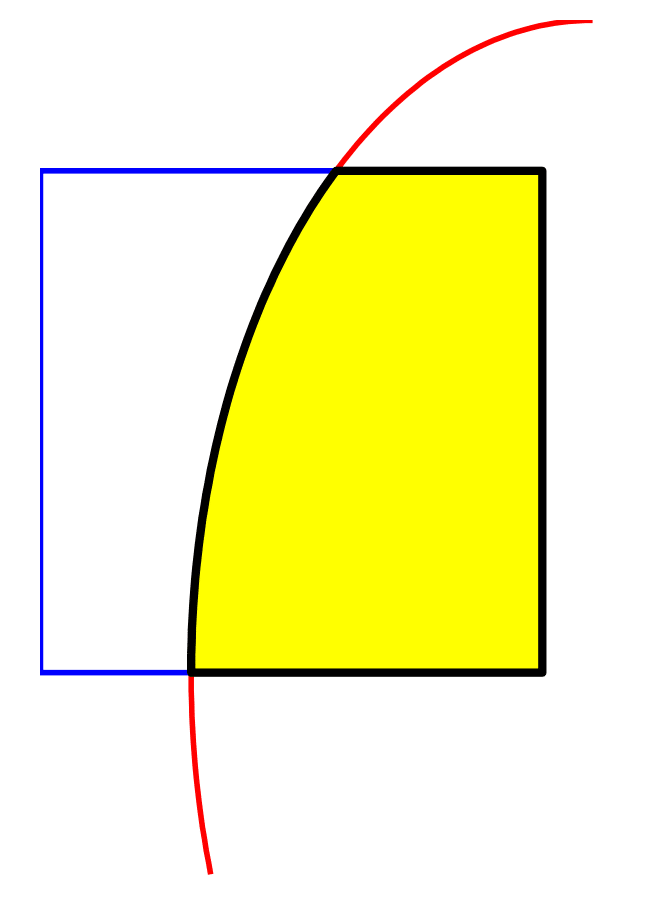} \\
   (c) & (d)
  \end{tabular}
\caption{\label{fig-eletypeB} Trimmed element of type B.}
\end{center}
\end{figure}

Except trimmed elements of type C, all trimmed elements are represented as quadrilaterals. For the curved quadrilateral which contains one curved edge as part of trimming curves, the mapping from unit rectangle is constructed as follows:
according to the location of the curved edge, there are four types of curved quadrilaterals as shown in Fig.~\ref{fig-eletypeB}. Suppose $u_2>u_1$, where $u_1,u_2$ are parameters of the trimming curve at the intersections. For each case, the mapping $Q$ between the curved quadrilateral and rectangle can be described as follows.

(a): if $u_1$ is the parameter of the left intersection point, the mapping $Q$ is constructed as
\begin{align}
\label{eq:1}
\begin{split}
X&=\phi_X\zeta+\frac{u-u_1}{u_2-u_1}(1-\zeta) \\
Y&=\phi_Y\zeta
\end{split}
\end{align}
otherwise,
\begin{align}
\tag{\ref{eq:1}$'$}
\begin{split}
X&=\phi_X\zeta+\frac{u_2-u}{u_2-u_1}(1-\zeta) \\
Y&=\phi_Y\zeta
\end{split}
\end{align}

(b): if $u_1$ is the  parameter of the right intersection point, the mapping $Q$ is constructed as
\begin{align}
\label{eq:2}
\begin{split}
X&=\phi_X(1-\zeta)+\frac{u-u_1}{u_2-u_1}\zeta \\
Y&=\phi_Y(1-\zeta)+\zeta
\end{split}
\end{align}
otherwise,
\begin{align}
\tag{\ref{eq:2}$'$}
\begin{split}
X&=\phi_X(1-\zeta)+\frac{u_2-u}{u_2-u_1}\zeta \\
Y&=\phi_Y(1-\zeta)+\zeta
\end{split}
\end{align}

(c): if $u_1$ is the parameter of the bottom intersection point, the mapping $Q$ is constructed as
\begin{align}
\label{eq:3}
\begin{split}
X&=\phi_X\zeta \\
Y&=\phi_Y\zeta+\frac{u-u_1}{u_2-u_1}(1-\zeta)
\end{split}
\end{align}
otherwise,
\begin{align}
\tag{\ref{eq:3}$'$}
\begin{split}
X&=\phi_X\zeta\\
Y&=\phi_Y\zeta+\frac{u_2-u}{u_2-u_1}(1-\zeta)
\end{split}
\end{align}

(d): if $u_1$ is the parameter of the top intersection point, the mapping $Q$ is constructed as
\begin{align}
\label{eq:4}
\begin{split}
X&=\phi_X(1-\zeta)+\zeta \\
Y&=\phi_Y(1-\zeta)+\frac{u-u_1}{u_2-u_1}\zeta
\end{split}
\end{align}
otherwise,
\begin{align}
\tag{\ref{eq:4}$'$}
\begin{split}
X&=\phi_X(1-\zeta)+\zeta\\
Y&=\phi_Y(1-\zeta)+\frac{u_2-u}{u_2-u_1}\zeta
\end{split}
\end{align}

\begin{figure}[!htb]
\small
\begin{center}
  \begin{tabular}{cc}
  \includegraphics[width=1.35in,clip]{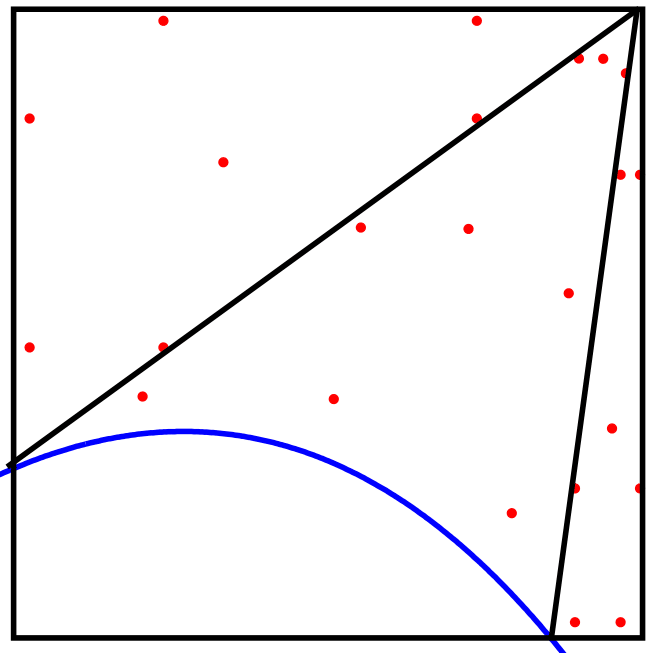} &
  \includegraphics[width=1.35in,clip]{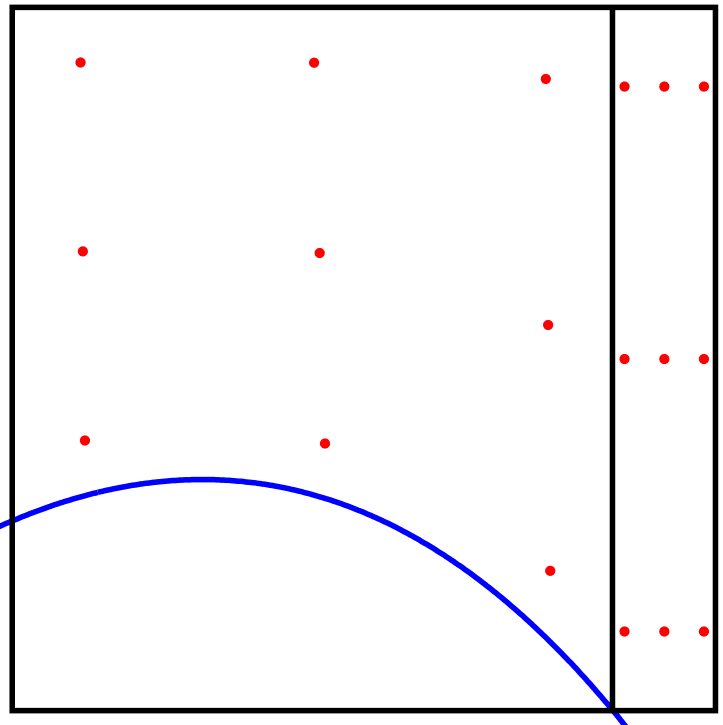} \\
  (a) & (b)\\
   \includegraphics[width=1.38in,clip]{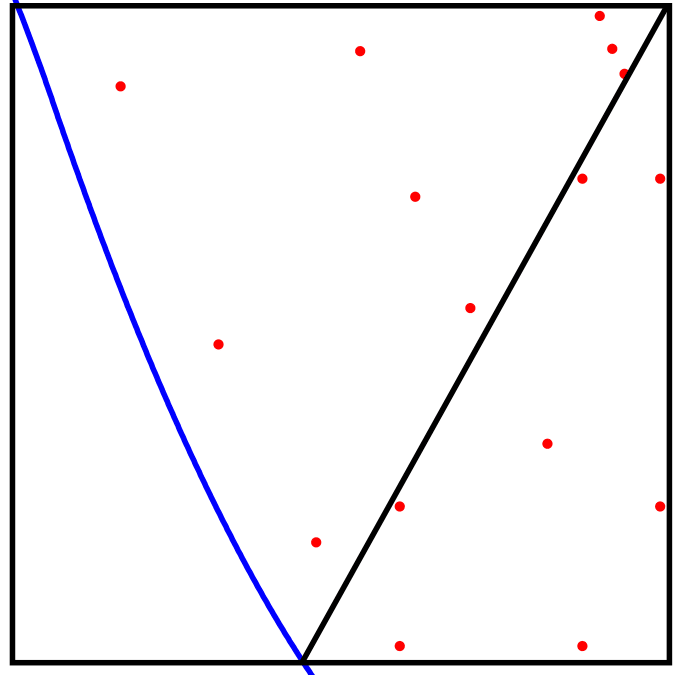} &
  \includegraphics[width=1.38in,clip]{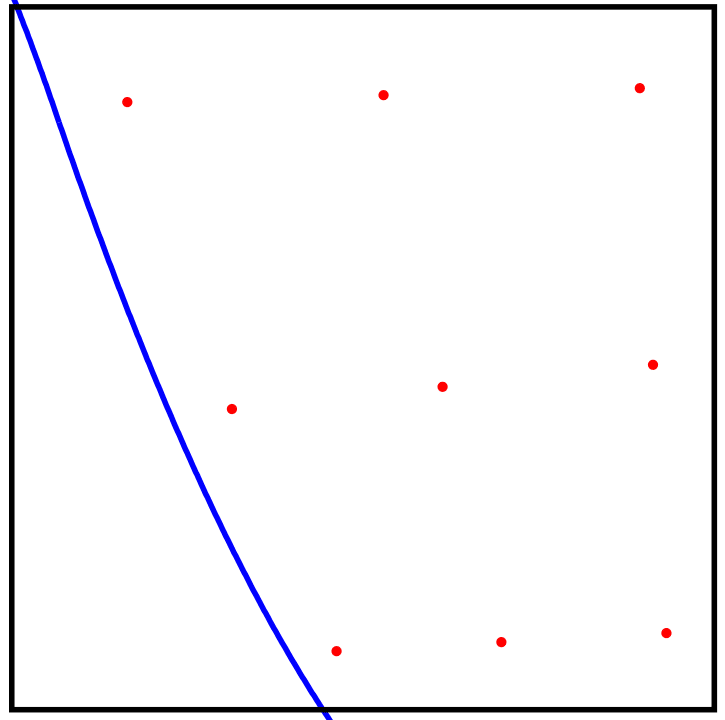} \\
  (c) & (d)
  \end{tabular}
\caption{\label{fig-compare-points} Integration points and segmentation of trimmed elements.(a)(c) type A and type B elements in ~\cite{Hyun-Jung Kim_2010}, (b)(d) type A and type B elements in our method.}
\end{center}
\end{figure}

Guass quadrature is commonly used in isogeometric FE approaches.
Compared to the method proposed by Kim et al.~\cite{Hyun-Jung Kim_2010}, the proposed method leads to less integration points. Fig.~\ref{fig-compare-points} shows the distribution of Gauss points in our approach compared to the approach in~\cite{Hyun-Jung Kim_2010} for one trimmed element. In fact, the reduction of integral points can be estimated. Suppose $n$ Gauss points are chosen for the normal triangle, and $m$ Gauss points are chosen for the curved triangle. As the number of integral points for quadrilateral element is the same with curved triangle, we can give the number of integral points for each type of trimmed element, see Table.~\ref{comparison-NoGPts}.
\begin{table}[!htb]
\centering
\caption{\label{comparison-NoGPts}Comparison of the number of integral points. }
\begin{tabular}{|c|c|c|}
  \hline
  element & original method & our method \\ \hline
  type A & 2n+m & 2m \\  \hline
  type B & n+m & m \\  \hline
  type C & m & m \\
  \hline
\end{tabular}
\end{table}

\section{Numerical examples}
\label{example}
In this section, we solve the Poisson equation on several trimmed geometries to show the effectiveness of our method, and  compare our results with results obtained by the method in~\cite{Hyun-Jung Kim_2010}.

The Poisson equation with Dirichlet boundary condition is given by,
\begin{equation}
\label{eq:poisson}
\left\{ \begin{aligned}
         -\Delta u &= f \\
       u|_{\partial \Omega}&=g,
      \end{aligned} \right.
\end{equation}
where $\Omega$ is the trimmed geometry represented by a NURBS surface and several trimming NURBS curves.

We use Lagrange multiplier method to impose the essential boundary condition as~\cite{Hyun-Jung Kim_2010} does. The Lagrange multipliers $\bm{\lambda}(u)$ are supposed to be expressed as
\[\bm{\lambda}(u)=\sum_{i=1}^lR_i^c(u)\lambda_i\]
where $R_i^c(u)$ is the NURBS basis function of trimming curve. The weak form with Lagrange multipliers is discretized as equations
\[\begin{aligned}
\bm{K}\bm{U}+\bm{A^T}\bm{\lambda}&=\bm{f} \\
\bm{A}\bm{U}&=\bm{b}
\end{aligned}\]
where
\[\begin{aligned}K_{ij}=\int_{\Omega}R_i^sR_j^sd\Omega, \quad f_i=\int_{\Omega}R_i^sfd\Omega\\
A_{ij}=\int_{\partial\Omega}R_i^cR_j^sd\partial\Omega, \quad b_i=\int_{\partial\Omega}R_i^cgd\partial\Omega
\end{aligned}\]
the $R_i^c(u)$ also denotes the NURBS basis function of trimming curve, $R_i^s(s,t)$ denotes the NURBS basis function of spline surface.
$\bigcup\Omega_i=\Omega, \Xi=\{\Omega_i\}_i$.
For subregion $\Omega_i$, we can use Gaussian points to compute the integration. More details can be found in~\cite{Hyun-Jung Kim_2010}.

For the method in ~\cite{Hyun-Jung Kim_2010}, basis functions are NURBS basis functions. Three Gaussian points are used in each direction for quadrilateral element, and seven Gaussian points are used for the regular triangle.

\begin{figure}[!htb]
\small
\begin{center}
  \begin{tabular}{ccc}
  \includegraphics[width=1.3in,clip]{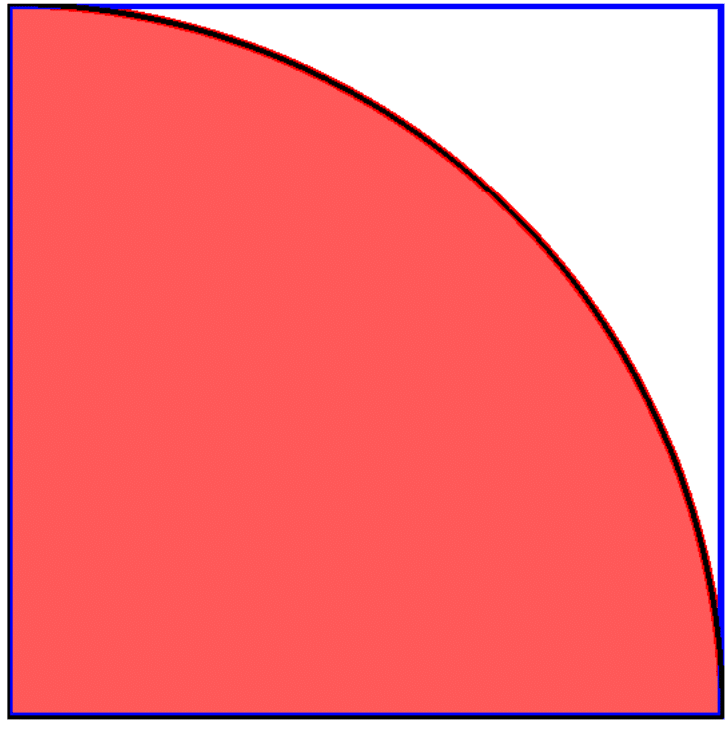} &
  \includegraphics[width=1.3in,clip]{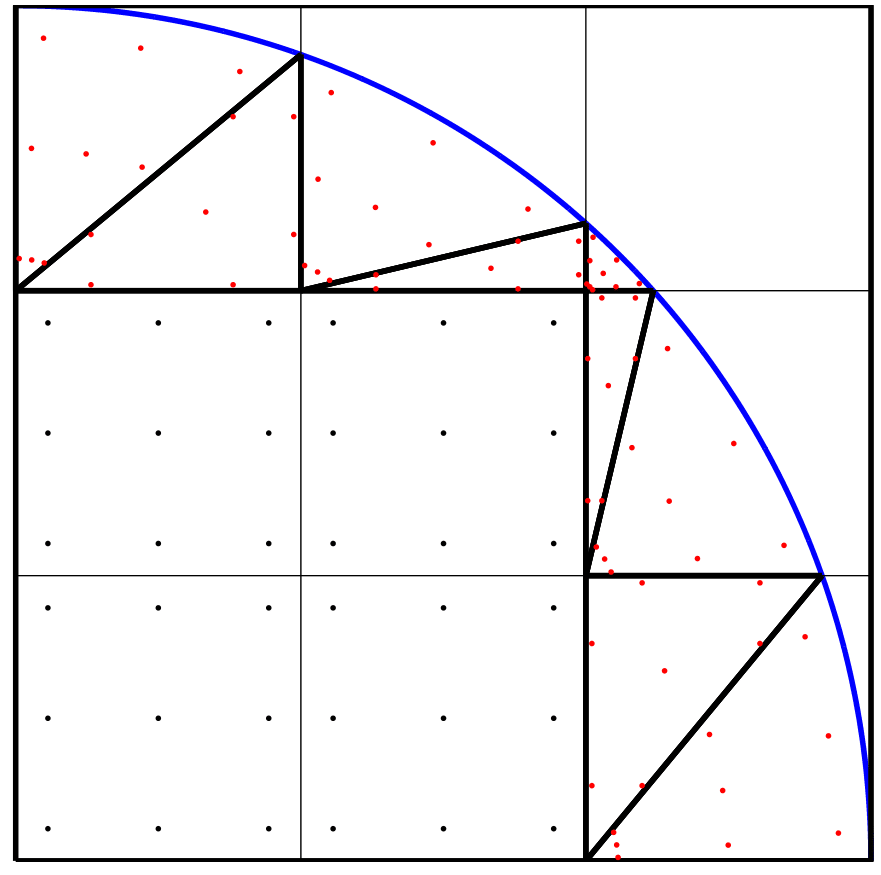} &
  \includegraphics[width=1.3in,clip]{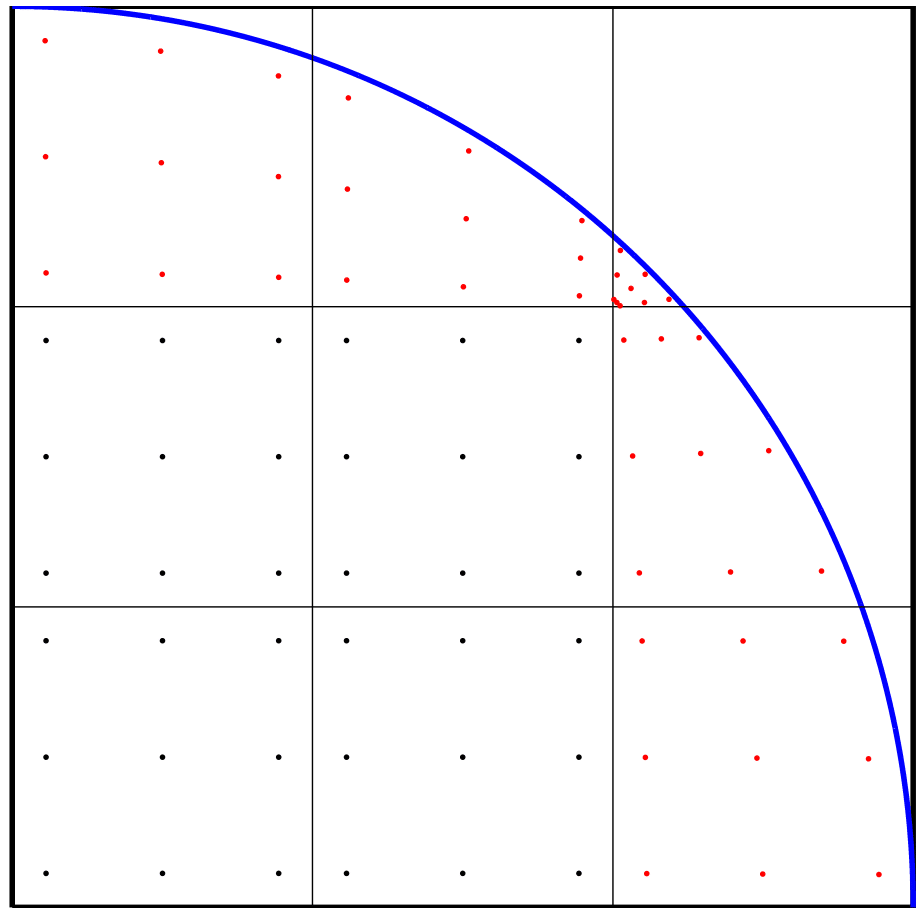} \\
  (a) & (b) & (c)
  \end{tabular}
\caption{\label{fig-ex1-intpts} Integration points and segmentation of trimmed elements with $3\times 3$ elements: (a) computational domain of EX1;(b) integration points in ~\cite{Hyun-Jung Kim_2010};(c) integration points in the proposed method.}
\end{center}
\end{figure}

\begin{figure}[!htb]
\small
\begin{center}
 \includegraphics[width=3in,clip]{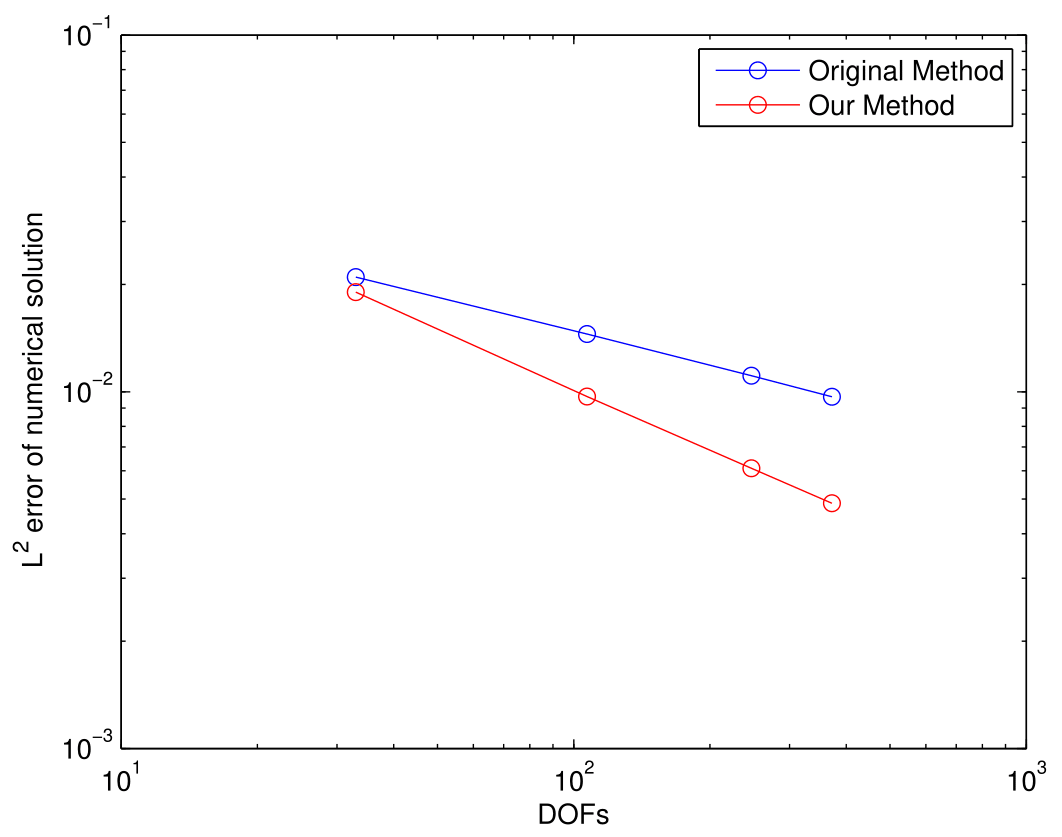}
\caption{\label{fig-errplot} Comparison of $L^2$ error.}
\end{center}
\end{figure}

In the first example, we choose a very simple geometry, the computational domain is fan-shaped. It is constructed by trimming a corner of a rectangle using an arc represented by a NURBS curve. In this example, we compare the corresponding $L^2$ error of numerical solution, and the condition number of stiffness matrix with the  method in ~\cite{Hyun-Jung Kim_2010} as shown in Table~\ref{comparison-ex1}. The error plot is shown in Fig.~\ref{fig-errplot} which is more illustrative. We also compare the computational cost of Ex1 as presented in Table~\ref{comparison-ex}, where $T_e$ represents the trimmed element and $\tilde{T_e}$ represents the integration element after the decomposition of $T_e$.

\begin{table}[htbp]
\small
\centering
\caption{\label{comparison-ex1}Comparison of our method with the proposed method in ~\cite{Hyun-Jung Kim_2010} for Ex1. }
\begin{tabular}{|c||c|c|c|c|}
\hline
\multirow{2}{*}{} &
\multicolumn{2}{|c|}{Method in ~\cite{Hyun-Jung Kim_2010}} & \multicolumn{2}{|c|}{Our method} \\
\cline{2-5}
 Number of element  & Cond. & $L^2$ error & Cond. & $L^2$ error \\
\hline
\hline
$5\times 5$  &  $7.5051\times 10^4$  & 0.0209621 &  $5.0469\times 10^4$ & 0.0190196 \\
\hline
$10\times 10$  &  $2.2417\times 10^{10}$  & 0.0145103 &  $1.5930\times 10^{10}$ & 0.00969233 \\
\hline
$20\times 20$  &  $4.3386\times 10^9$  & 0.0096861 &  $2.9344\times 10^9$ & 0.00487051 \\
\hline
\end{tabular}
\end{table}

\begin{figure}[!htb]
\small
\begin{center}
  \begin{tabular}{ccc}
  \includegraphics[width=1.4in,clip]{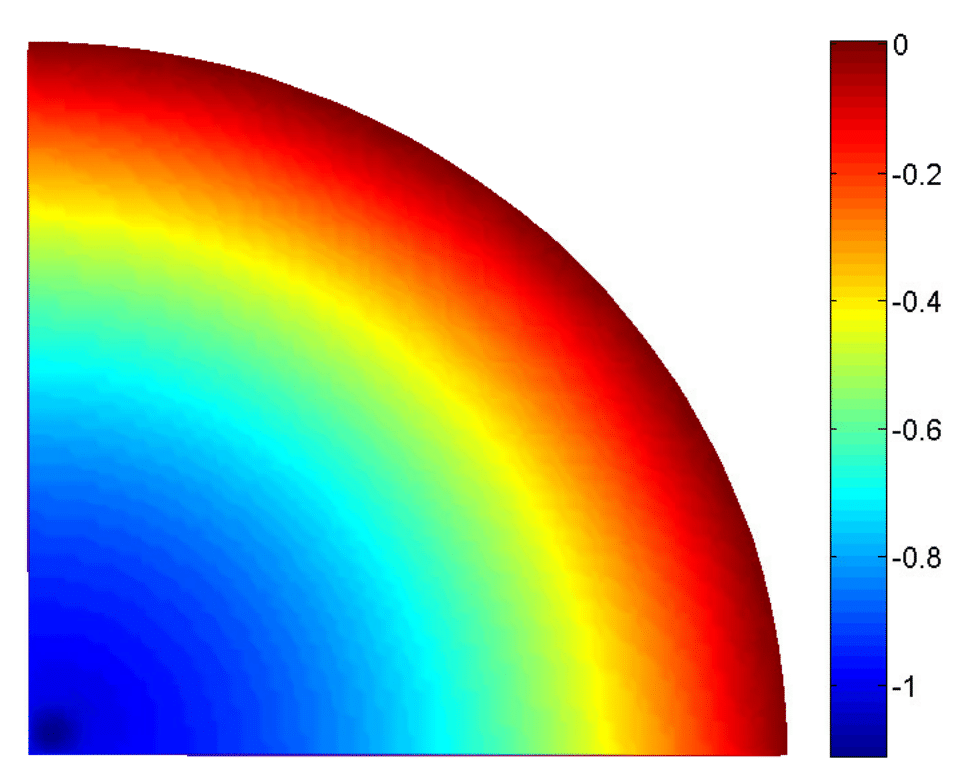} &
  \includegraphics[width=1.4in,clip]{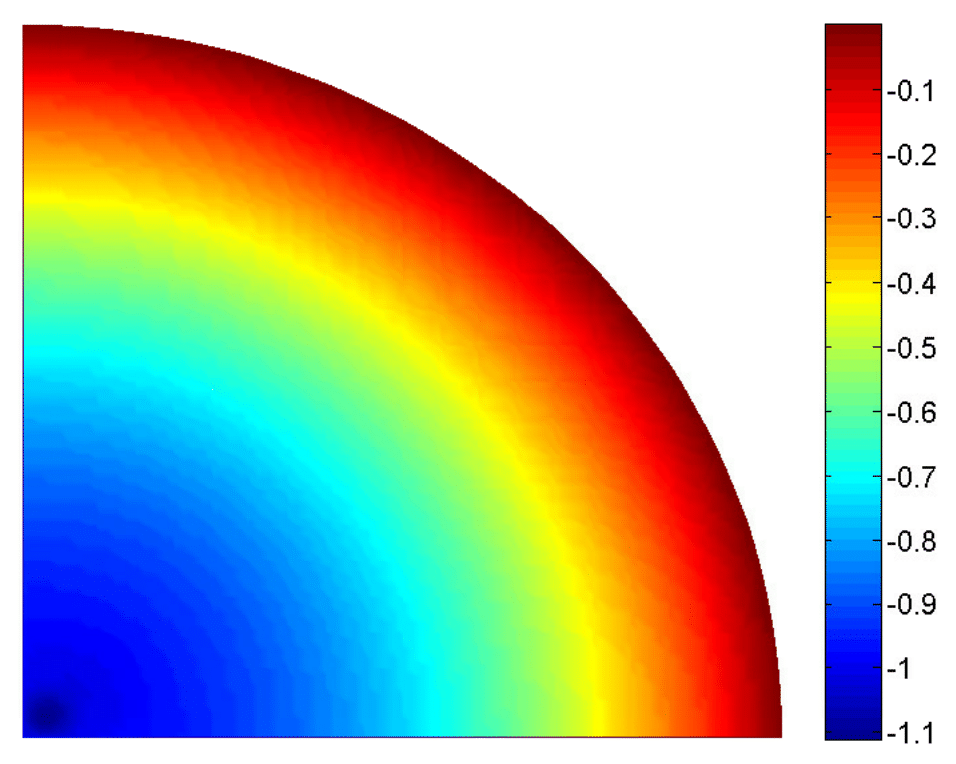} &
  \includegraphics[width=1.4in,clip]{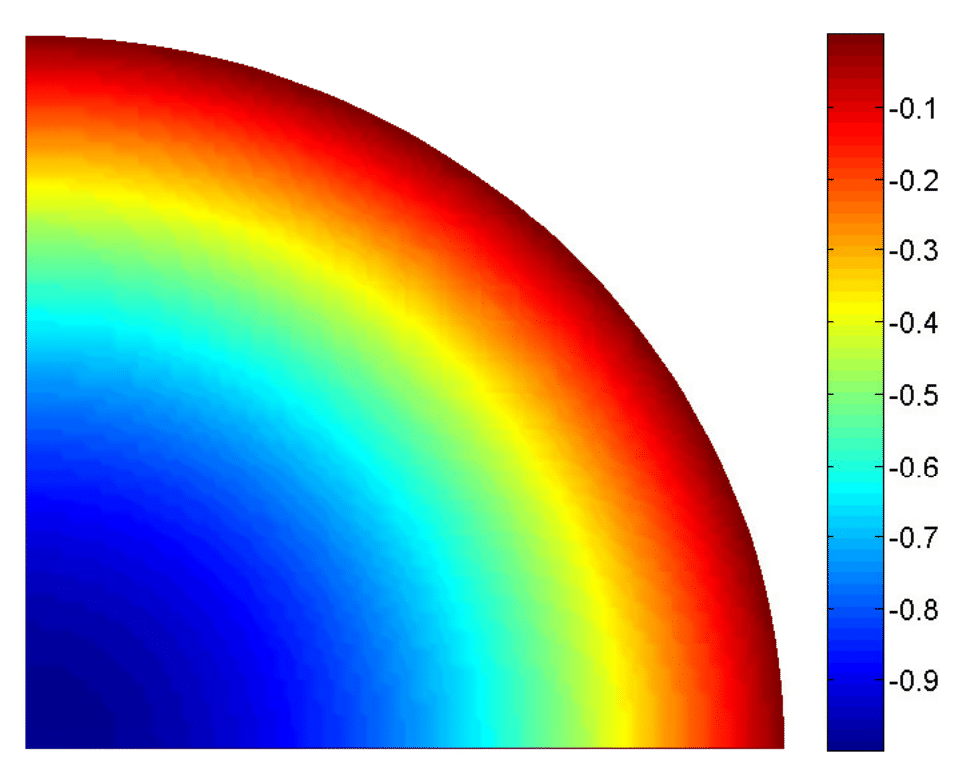} \\
  (a) & (b) & (c)
  \end{tabular}
\caption{\label{fig-ex1-result} Comparison of numerical solution  with $20\times 20$ grid.(a) the  method in ~\cite{Hyun-Jung Kim_2010};(b) our method;(c)exact solution.}
\end{center}
\end{figure}

We use the method presented in \cite{Hyun-Jung Kim_2009} to find all the active elements, and construct the mapping from the unit square $[0,1]\times [0,1] $ to each trimmed element as described in section~\ref{integration-pts}. When the spline surface consists of $3\times 3$ elements, with our method, the distribution of the integration  points on the computational domain is more regular  than the method in ~\cite{Hyun-Jung Kim_2010} as illustrated in Fig.~\ref{fig-ex1-intpts}. The corresponding numerical solution is shown in Fig.~\ref{fig-ex1-result}. From Table~\ref{comparison-ex1}, we can see that the condition number of the stiffness matrix and the $L^2$ error of the numerical solution is reduced almost by one half compared to  the  method in  ~\cite{Hyun-Jung Kim_2010} on the refined grid.

 \begin{figure}[!htb]
\small
\begin{center}
  \begin{tabular}{cc}
  \includegraphics[width=1.4in,clip]{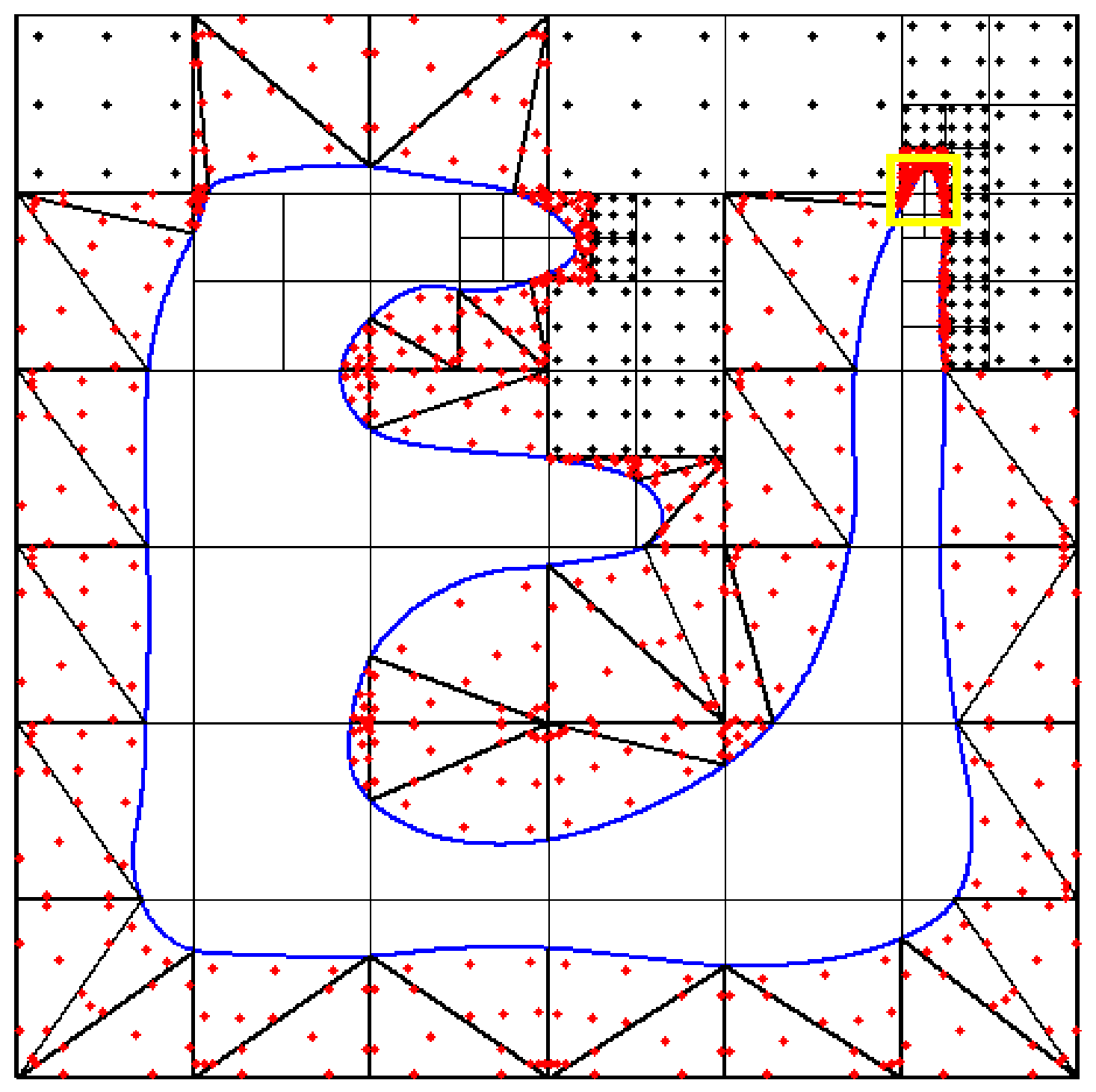} &
  \includegraphics[width=1.4in,clip]{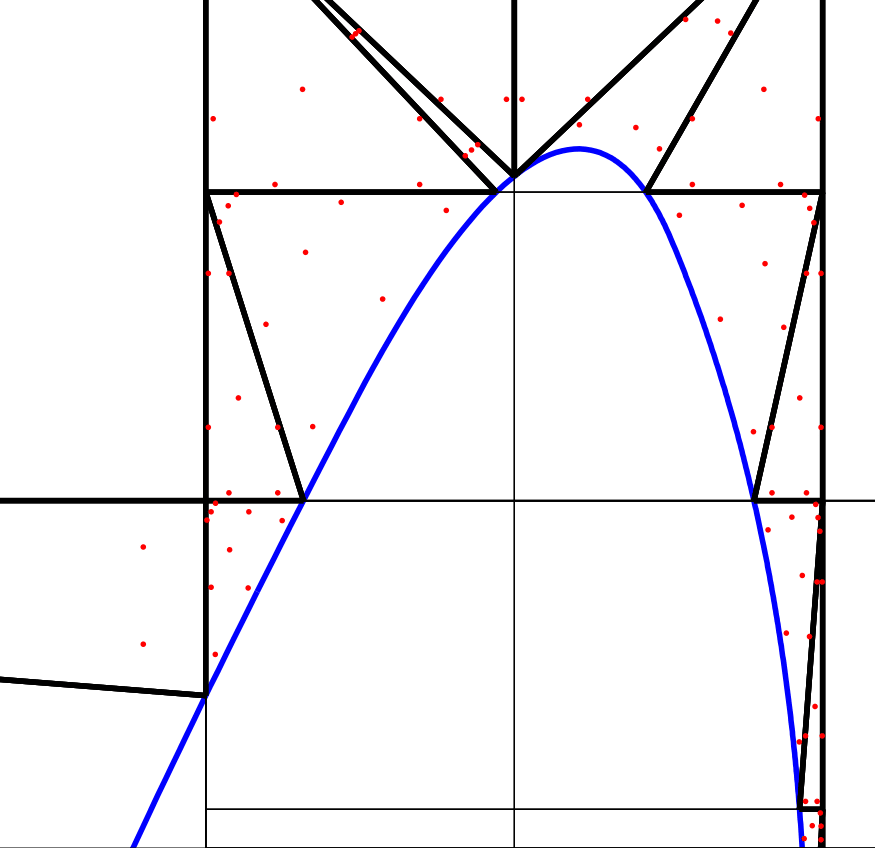} \\
  (a) & (b)
  \end{tabular}
  \begin{tabular}{cc}
  \includegraphics[width=1.4in,clip]{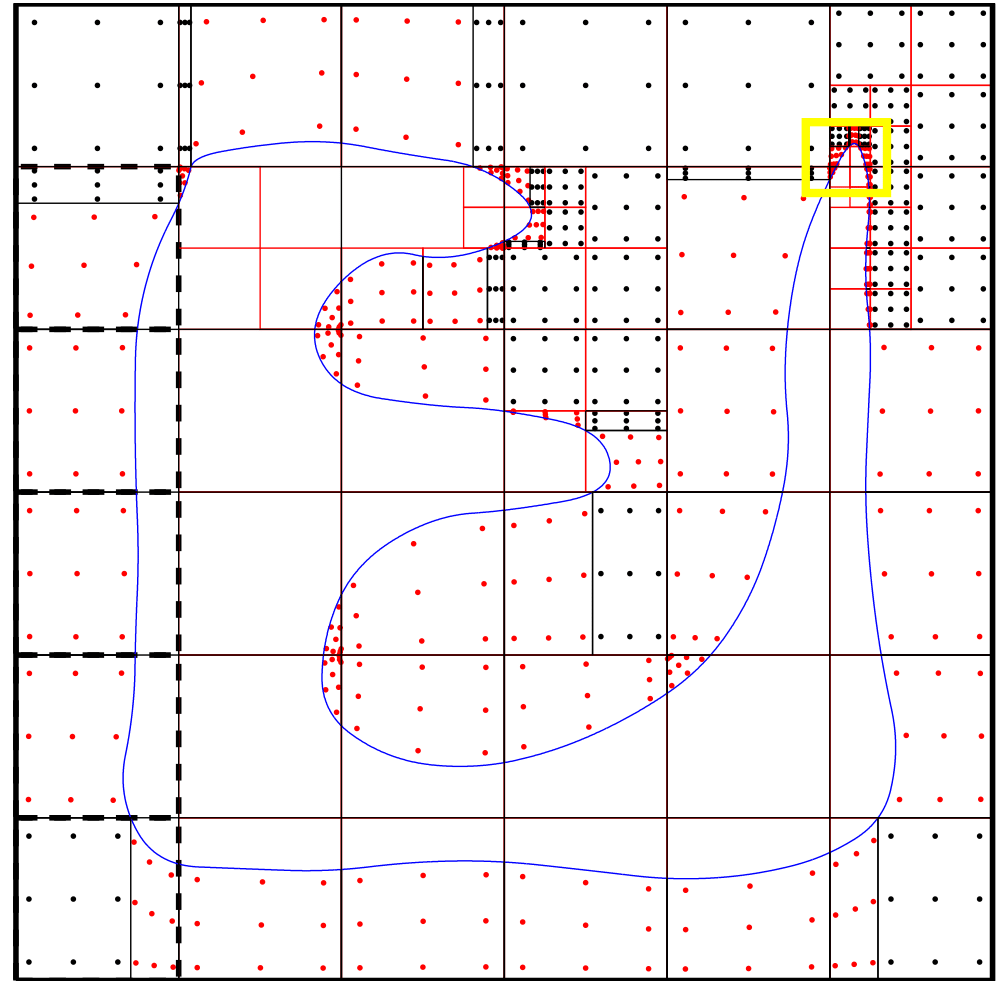} &
  \includegraphics[width=1.4in,clip]{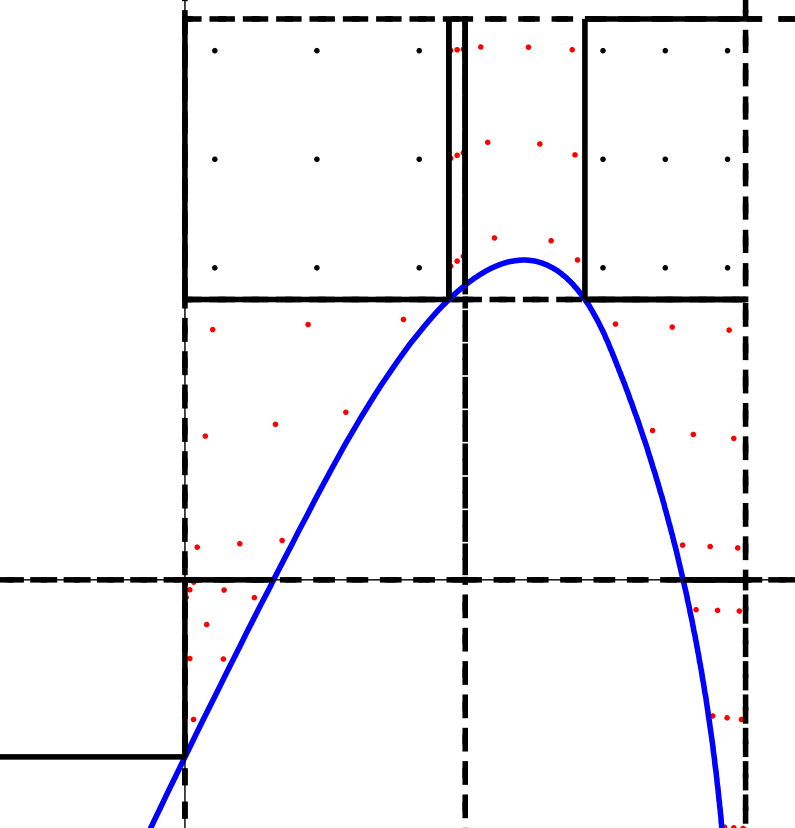} \\
  (c) & (d)
  \end{tabular}
\caption{\label{fig-ex2-intpts} The computational domain of EX2. (a) elements and integral points of the  method in ~\cite{Hyun-Jung Kim_2010};(b)enlarge the area of yellow rectangle in (a);(c) elements and integral points of our method;(d)enlarge the area of yellow rectangle in (c).}
\end{center}
\end{figure}

 \begin{figure}[!htb]
\small
\begin{center}
  \begin{tabular}{ccc}
  \includegraphics[width=1.4in,clip]{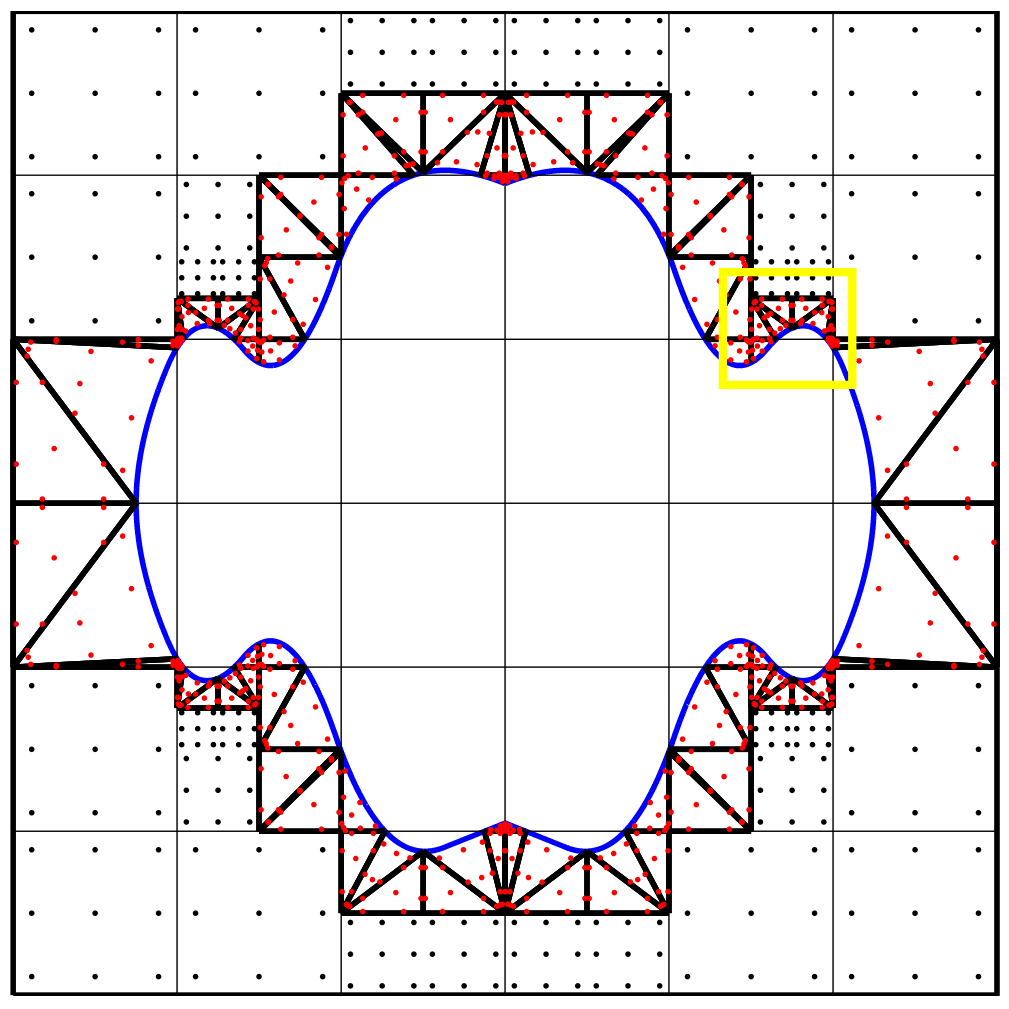} &
  \includegraphics[width=1.4in,clip]{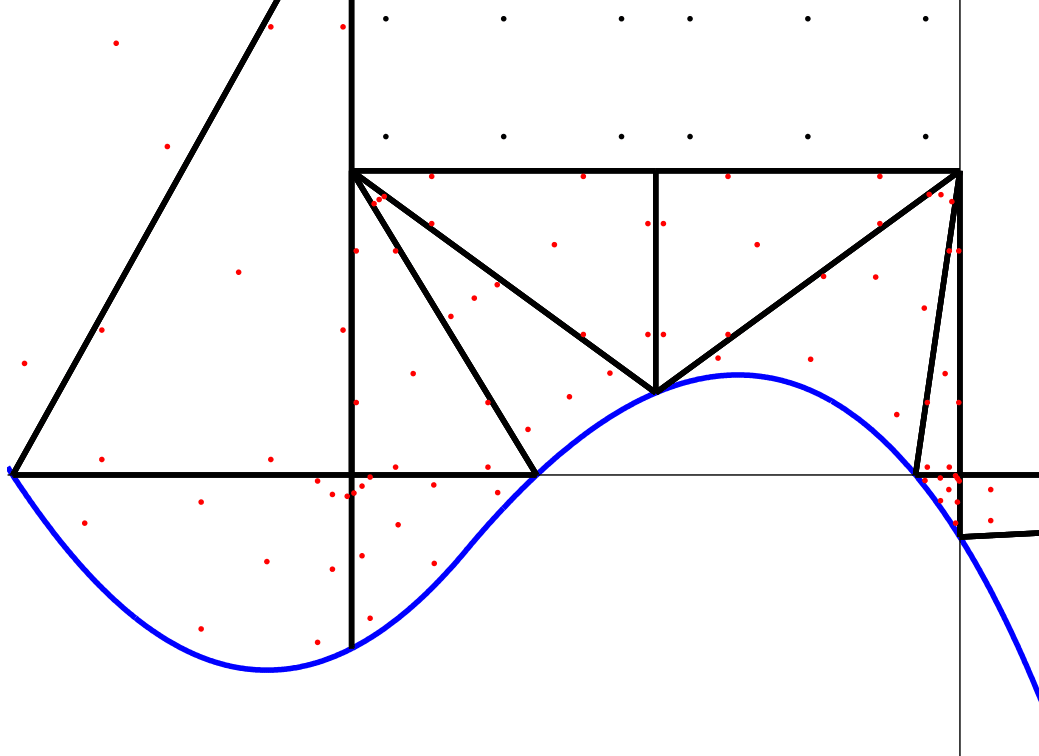} \\
  (a) & (b)
  \end{tabular}
  \begin{tabular}{cc}
  \includegraphics[width=1.4in,clip]{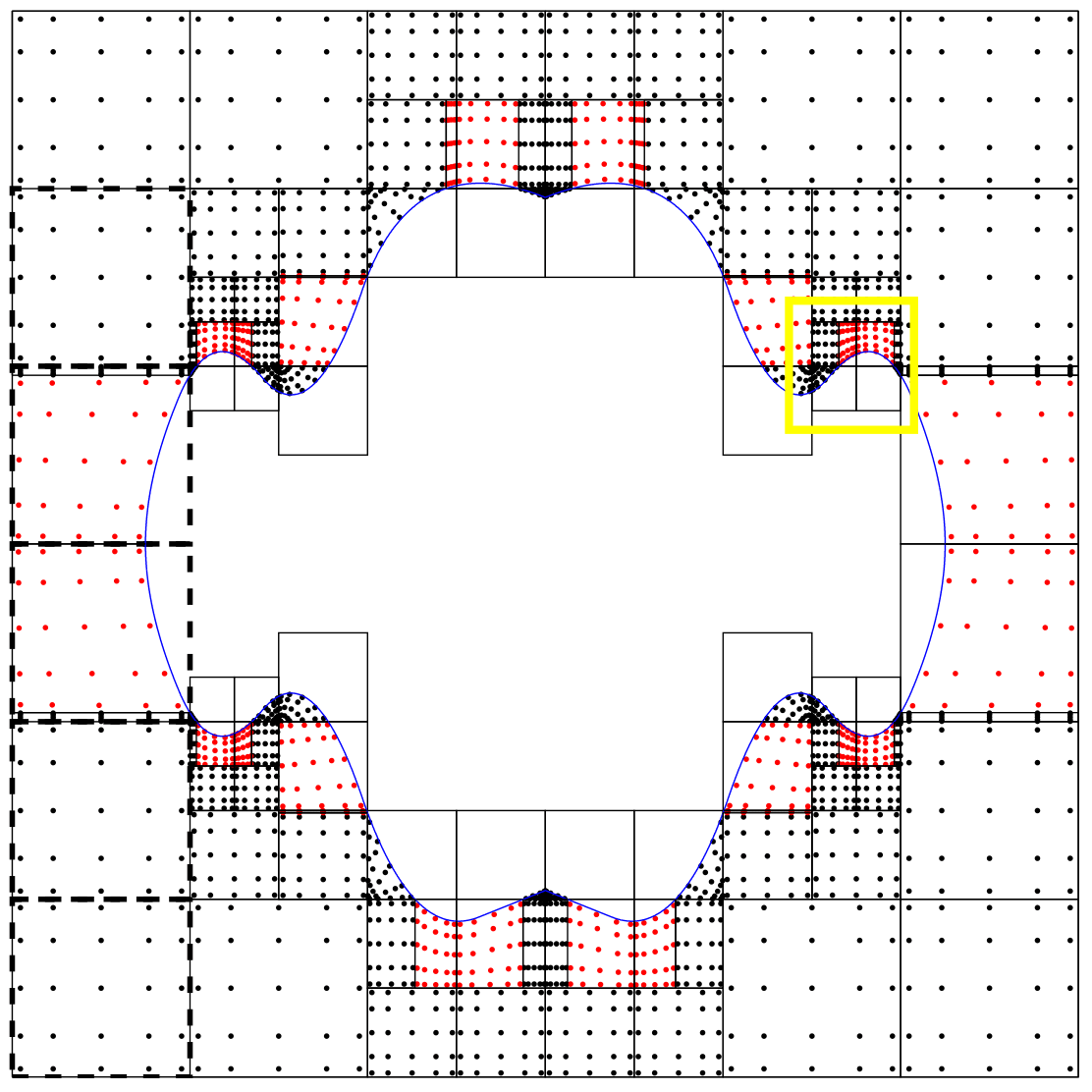} &
  \includegraphics[width=1.4in,clip]{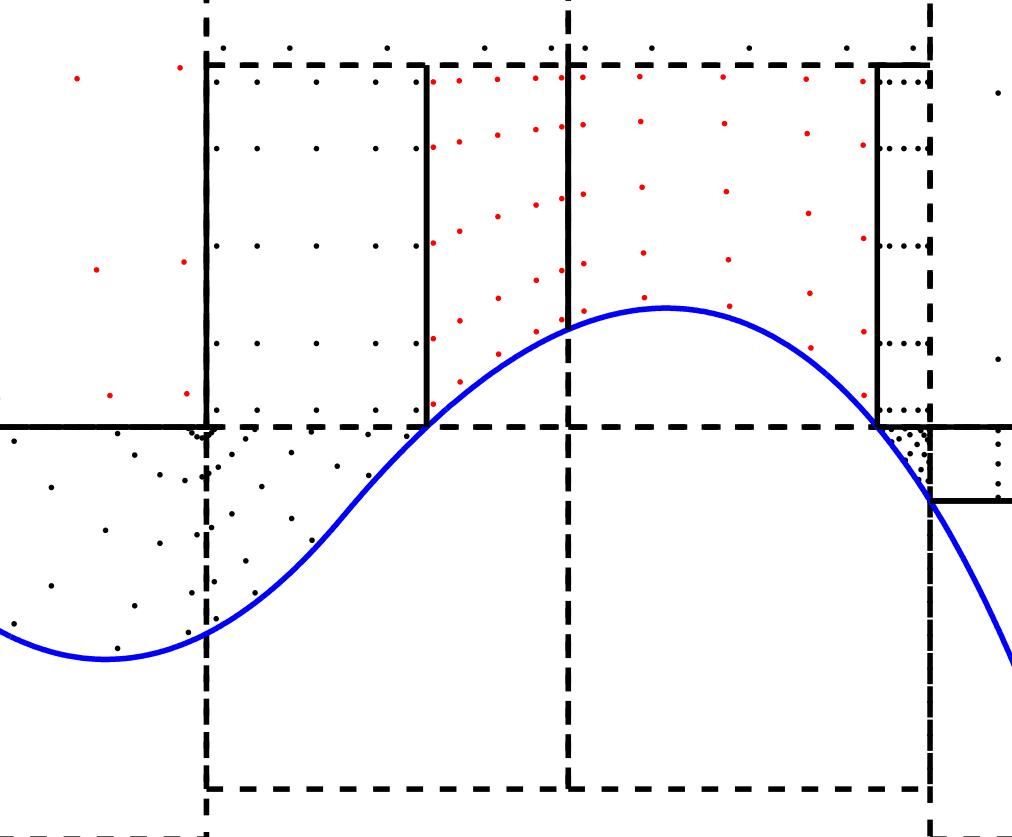} \\
  (c) & (d)
  \end{tabular}
    \begin{tabular}{cc}
  \includegraphics[width=1.5in,clip]{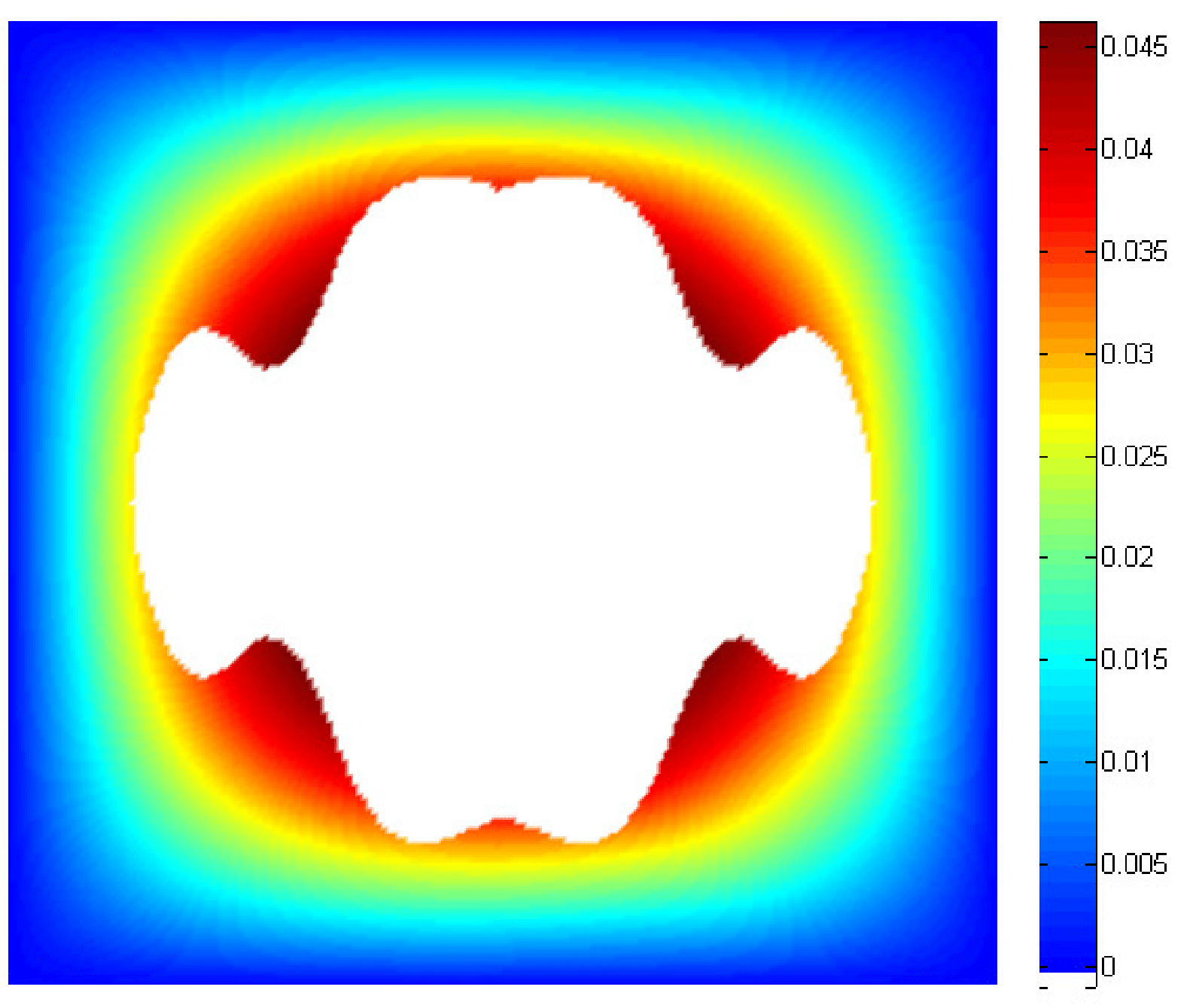} &
  \includegraphics[width=1.5in,clip]{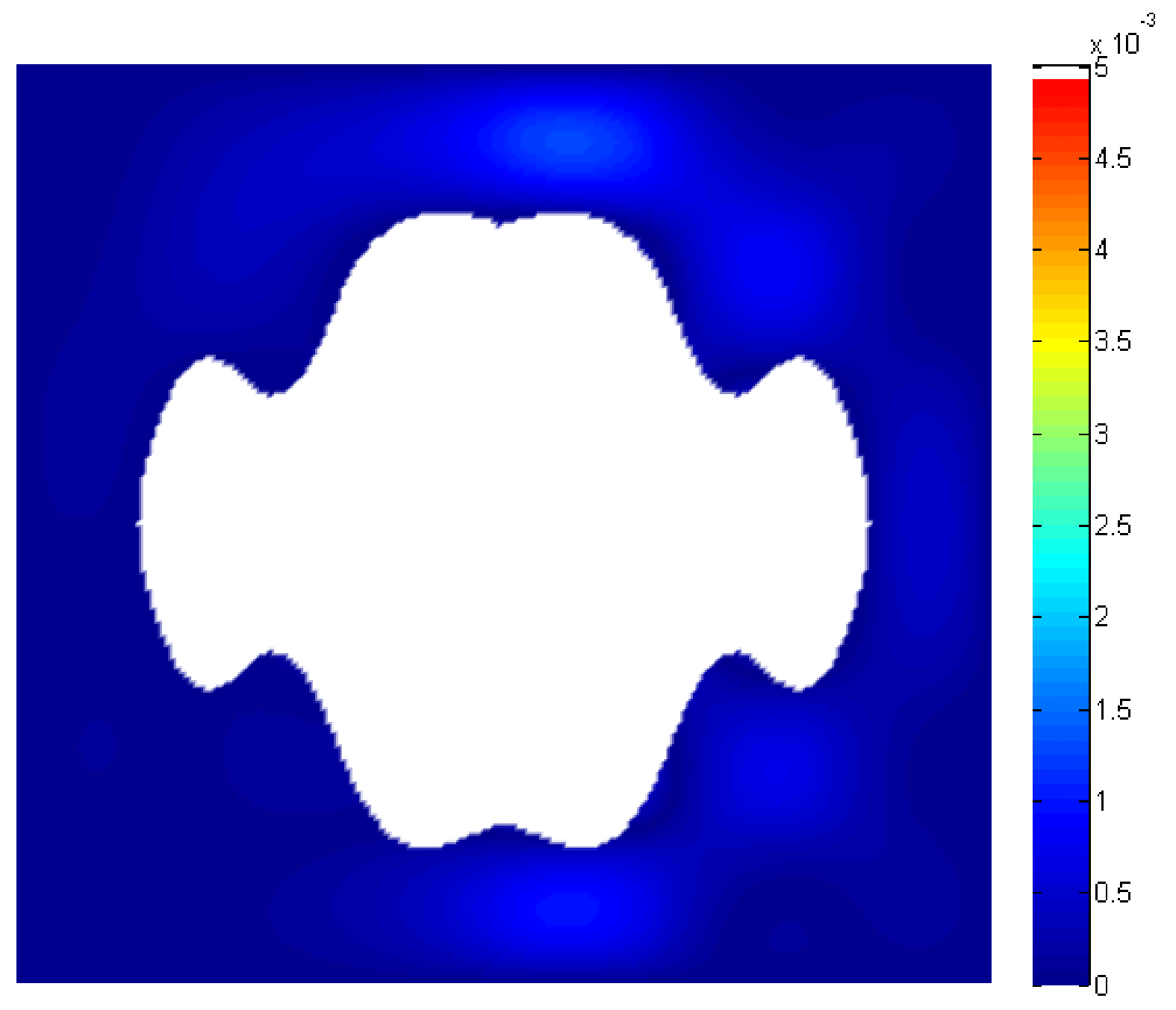} \\
  (e) & (f)
  \end{tabular}
\caption{\label{fig-ex3-intpts} The computational domain of EX3: (a)  elements and integral points of the  method in ~\cite{Hyun-Jung Kim_2010};(b)enlarge the area of yellow rectangle in (a);(c) elements and integral points of our method;(d) enlarge the area of yellow rectangle in (c);(e)numerical solution;(f)$L^2$ error is $4.43065\times 10^{-4}$. }
\end{center}
\end{figure}

In the second example, we construct a computational domain with complex geometry, where the rectangle is trimmed by a closed spline curve. There are two protrusions in the interior of the final trimmed geometry. In this example, the surface contains $6\times 6$ elements first. However, there are other kind of trimmed elements except of three types we processed in this coarse mesh, so local refinement is performed on the surface as described in~\cite{Hyun-Jung Kim_2009} until there are only three types of trimmed elements. In this example, there are many  type B trimmed elements, as shown in Fig.~\ref{fig-ex2-intpts}  . The element of type B is decomposed into  two triangles with the method in ~\cite{Hyun-Jung Kim_2010}, integration  on this element then becomes integration on these two triangles. In the proposed method, we construct a mapping from type B element to rectangle directly while keeping the number of integral elements. Our  method  can reduce a half integral points and integral elements for this type of trimmed element. For type A,  one third of integral points and integral elements can be reduced by our method.

In the third example, the number of trimmed elements of type A and type C are more than type B, in this case the reduction of integral points and integral elements is not as significant as the first example. It can be clearly observed from Table 2.

A round hole is trimmed out from a rectangle as computational domain in the last example. In this example, the number of trimmed elements of type A becomes more and more during mesh refinement process , and the reduction of integral elements is clearly demonstrated.

 \begin{figure}[!htb]
\small
\begin{center}
  \begin{tabular}{cc}
  \includegraphics[width=1.4in,clip]{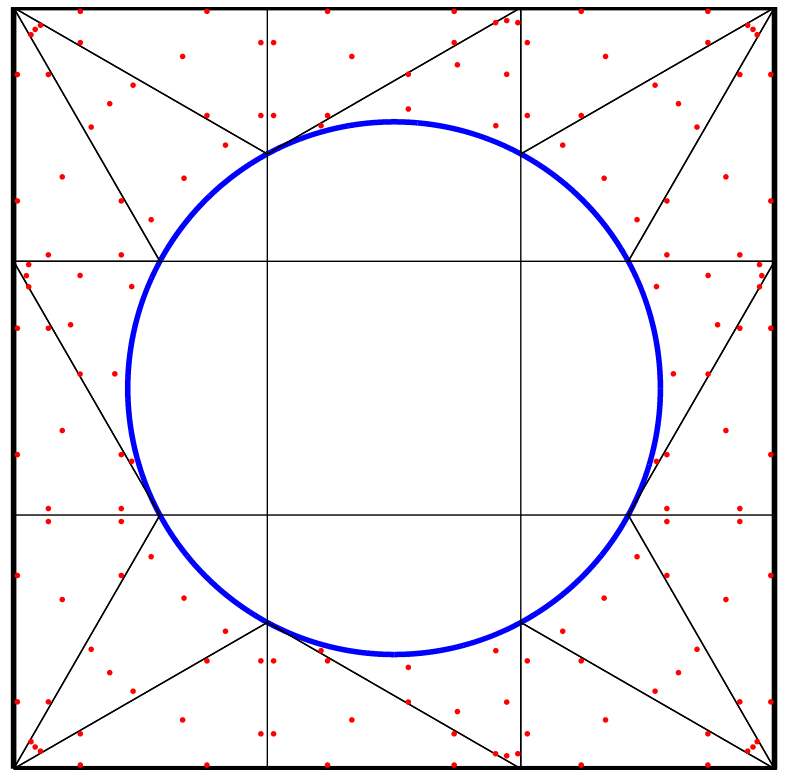} &
  \includegraphics[width=1.4in,clip]{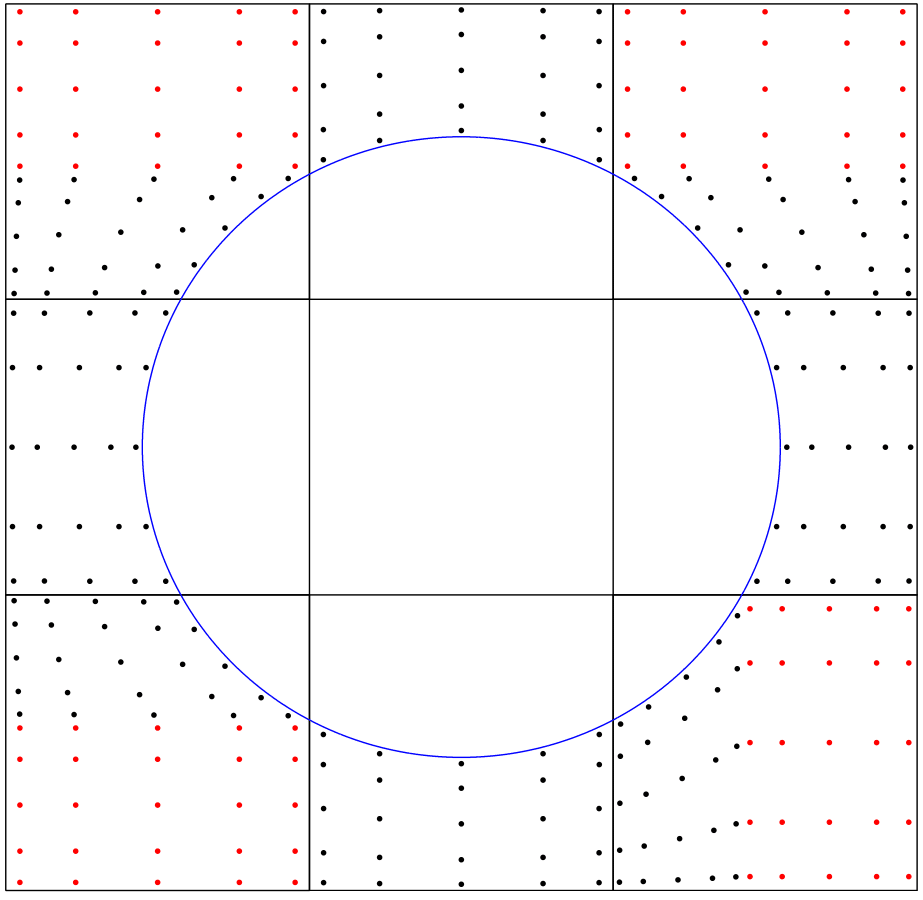} \\
  (a) & (b)
  \end{tabular}
   \begin{tabular}{cc}
  \includegraphics[width=1.5in,clip]{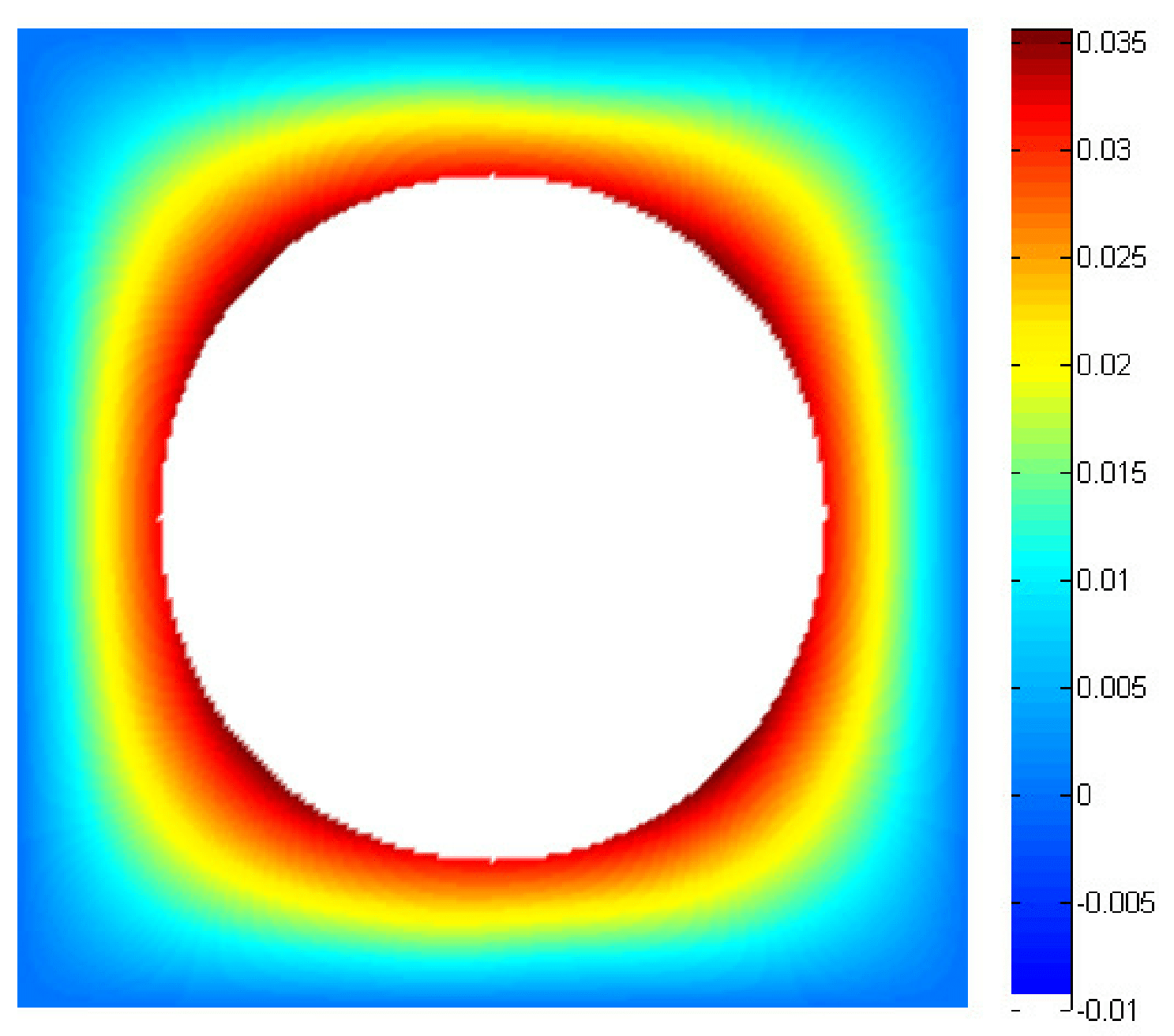} &
  \includegraphics[width=1.5in,clip]{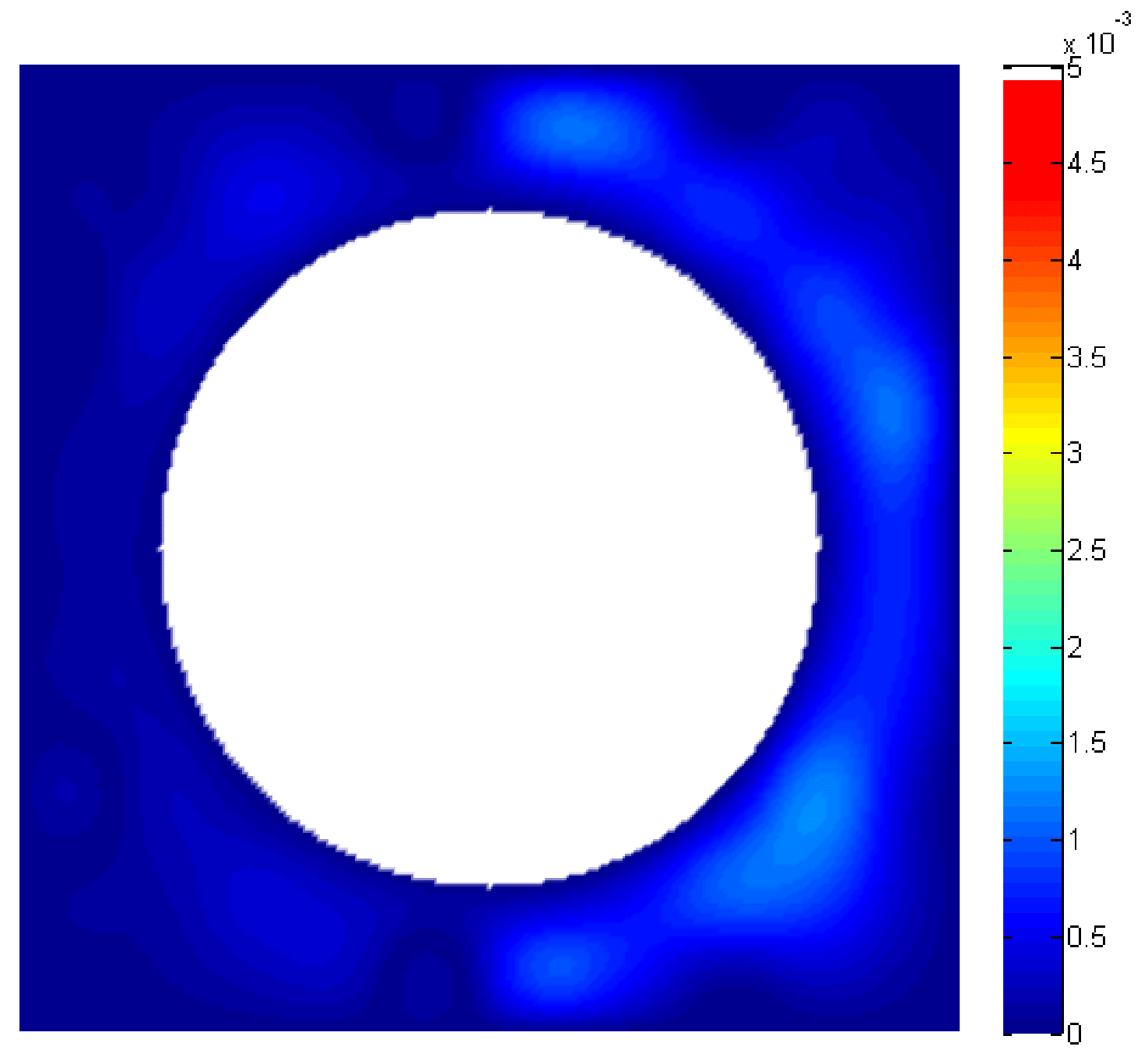} \\
  (c) & (d)
  \end{tabular}
\caption{\label{fig-ex2-intpts}  The computational domain of EX4. (a)elements and integral points of the method in ~\cite{Hyun-Jung Kim_2010};(b)elements and integral points of our method; (c)numerical solution of $10\times 10$ grid;(d)$L^2$ error is $3.26665\times 10^{-4}$.}
\end{center}
\end{figure}

 \begin{table}[htbp]
\centering
\caption{\label{comparison-ex}Comparison of computational cost with the Method in ~\cite{Hyun-Jung Kim_2010}}
\small
\begin{tabular}{|c||c|c|c|c|c|c|}
\hline
\multirow{2}{*}{} & \multirow{2}{*}{Mesh Size} & \multirow{2}{*}{Number of $T_e$} &
\multicolumn{2}{|c|}{Number of $\tilde{T}_e$} & \multicolumn{2}{|c|}{Number of integral points in $T_e$} \\
\cline{4-7}
& & & \small{Method in ~\cite{Hyun-Jung Kim_2010}} & \small{Our method} &  \small{Method in ~\cite{Hyun-Jung Kim_2010}} & \small{Our method}\\
    \hline 
    \hline
   Ex1  & $3\times 3$&5   & 9  & 5 & 73 & 45 \\
    & $10\times 10$&17   & 31  & 19 & 251 & 171 \\
    & $20\times 20$&37   & 71  & 45 & 571 & 405 \\
\hline
Ex2  & $6\times 6$  & 47 & 98 & 60 & 780 & 486 \\
\hline
Ex3  & $6\times 6$  & 48 & 100 & 72 & 796 & 648 \\
\hline
  Ex4   & $3\times 3$&8   & 20  & 12 & 156 & 108 \\
    & $10\times 10$&28   & 60  & 40 & 476 & 360 \\
    & $20\times 20$&56   & 112  & 76 & 896 & 684 \\
\hline
\end{tabular}
\end{table}

\section{Conclusion}
\label{conlude}

In this paper, we  propose an improved method of isogeometric analysis over trimmed geometries  on  two-dimensional planar computational domain. By the proposed method, the  integral elements and integral points in analysis process can be reduced significantly, which improves the efficiency of analysis. Moreover, compared with the previous method, the distribution of integral points is more regular, and the accuracy of numerical solution is also improved. Extension to three-dimensional isogeometric analysis ~\cite{xu:cad17,xu:spm2012} will be a part of future work.

\section*{Acknowledgment}

This research was supported by the National Nature Science Foundation of China under Grant Nos. 61602138 and 61472111, Zhejiang Provincial Natural Science Foundation of China under Grant Nos. LQ16F020005 and LR16F020003, and the Open Project Program of the State Key Lab of CAD\&CG (A1703), Zhejiang University.

\section*{References}

\bibliographystyle{abbv}

\end{document}